


\magnification\magstep1
\baselineskip14.75pt
\advance\vsize1pt
\newread\AUX\immediate\openin\AUX=\jobname.aux
\newcount\relFnno
\def\ref#1{\expandafter\edef\csname#1\endcsname}
\ifeof\AUX\immediate\write16{\jobname.aux gibt es nicht!}\else
\input \jobname.aux
\fi\immediate\closein\AUX



\def\ignore{\bgroup
\catcode`\;=0\catcode`\^^I=14\catcode`\^^J=14\catcode`\^^M=14
\catcode`\ =14\catcode`\!=14\catcode`\"=14\catcode`\#=14\catcode`\$=14
\catcode`\&=14\catcode`\'=14\catcode`\(=14\catcode`\)=14\catcode`\*=14
\catcode`+=14\catcode`\,=14\catcode`\-=14\catcode`\.=14\catcode`\/=14
\catcode`\0=14\catcode`\1=14\catcode`\2=14\catcode`\3=14\catcode`\4=14
\catcode`\5=14\catcode`\6=14\catcode`\7=14\catcode`\8=14\catcode`\9=14
\catcode`\:=14\catcode`\<=14\catcode`\==14\catcode`\>=14\catcode`\?=14
\catcode`\@=14\catcode`\A=14\catcode`\B=14\catcode`\C=14\catcode`\D=14
\catcode`\E=14\catcode`\F=14\catcode`\G=14\catcode`\H=14\catcode`\I=14
\catcode`\J=14\catcode`\K=14\catcode`\L=14\catcode`\M=14\catcode`\N=14
\catcode`\O=14\catcode`\P=14\catcode`\Q=14\catcode`\R=14\catcode`\S=14
\catcode`\T=14\catcode`\U=14\catcode`\V=14\catcode`\W=14\catcode`\X=14
\catcode`\Y=14\catcode`\Z=14\catcode`\[=14\catcode`\\=14\catcode`\]=14
\catcode`\^=14\catcode`\_=14\catcode`\`=14\catcode`\a=14\catcode`\b=14
\catcode`\c=14\catcode`\d=14\catcode`\e=14\catcode`\f=14\catcode`\g=14
\catcode`\h=14\catcode`\i=14\catcode`\j=14\catcode`\k=14\catcode`\l=14
\catcode`\m=14\catcode`\n=14\catcode`\o=14\catcode`\p=14\catcode`\q=14
\catcode`\r=14\catcode`\s=14\catcode`\t=14\catcode`\u=14\catcode`\v=14
\catcode`\w=14\catcode`\x=14\catcode`\y=14\catcode`\z=14\catcode`\{=14
\catcode`\|=14\catcode`\}=14\catcode`\~=14\catcode`\^^?=14
\Ignoriere}
\def\Ignoriere#1\;{\egroup}

\newcount\itemcount
\def\resetitem{\global\itemcount0}\resetitem
\newcount\Itemcount
\Itemcount0
\newcount\maxItemcount
\maxItemcount=0

\def\FILTER\fam\itfam\tenit#1){#1}

\def\Item#1{\global\advance\itemcount1
\edef\TEXT{{\it\romannumeral\itemcount)}}%
\ifx?#1?\relax\else
\ifnum#1>\maxItemcount\global\maxItemcount=#1\fi
\expandafter\ifx\csname I#1\endcsname\relax\else
\edef\testA{\csname I#1\endcsname}
\expandafter\expandafter\def\expandafter\testB\testA
\edef\testC{\expandafter\FILTER\testB}
\edef\testD{\csname0\testC0\endcsname}\fi
\edef\testE{\csname0\romannumeral\itemcount0\endcsname}
\ifx\testD\testE\relax\else
\immediate\write16{I#1 hat sich geaendert!}\fi
\expandwrite\AUX{\neverexpand\ref{I#1}{\TEXT}}\fi
\item{\ifx?#1?\relax\else\marginnote{I#1}\fi\TEXT}}

\def\today{\number\day.~\ifcase\month\or
  Januar\or Februar\or M{\"a}rz\or April\or Mai\or Juni\or
  Juli\or August\or September\or Oktober\or November\or Dezember\fi
  \space\number\year}
\font\sevenex=cmex7
\scriptfont3=\sevenex
\font\fiveex=cmex10 scaled 500
\scriptscriptfont3=\fiveex
\def\A{{\bf A}}
\def\G{{\bf G}}
\def\P{{\bf P}}

\def\phi{\varphi}
\def\epsilon{\varepsilon}
\def\theta{\vartheta}
\def\uauf{\lower1.7pt\hbox to 3pt{%
\vbox{\offinterlineskip
\hbox{\vbox to 8.5pt{\leaders\vrule width0.2pt\vfill}%
\kern-.3pt\hbox{\lams\char"76}\kern-0.3pt%
$\raise1pt\hbox{\lams\char"76}$}}\hfil}}

\def\title#1{\par
{\baselineskip1.5\baselineskip\rightskip0pt plus 5truecm
\leavevmode\vskip0truecm\noindent\font\BF=cmbx10 scaled \magstep2\BF #1\par}
\vskip1truecm
\leftline{\font\CSC=cmcsc10
{\CSC Friedrich Knop}}
\leftline{Department of Mathematics, Rutgers University, New Brunswick NJ
08903, USA}
\leftline{knop@math.rutgers.edu}
\vskip1truecm
\par}

\def\cite#1{\expandafter\ifx\csname#1\endcsname\relax
{\bf?}\immediate\write16{#1 ist nicht definiert!}\else\csname#1\endcsname\fi}
\def\expandwrite#1#2{\edef\next{\write#1{#2}}\next}
\def\neverexpand{\noexpand\noexpand\noexpand}
\def\strip#1\ {}
\def\ncite#1{\expandafter\ifx\csname#1\endcsname\relax
{\bf?}\immediate\write16{#1 ist nicht definiert!}\else
\expandafter\expandafter\expandafter\strip\csname#1\endcsname\fi}
\newwrite\AUX
\immediate\openout\AUX=\jobname.aux
\font\eightrm=cmr8\font\sixrm=cmr6
\font\eighti=cmmi8
\font\eightit=cmti8
\font\eightbf=cmbx8
\font\eightcsc=cmcsc10 scaled 833
\def\eightpoint{%
\textfont0=\eightrm\scriptfont0=\sixrm\def\rm{\fam0\eightrm}%
\textfont1=\eighti
\textfont\bffam=\eightbf\def\bf{\fam\bffam\eightbf}%
\textfont\itfam=\eightit\def\it{\fam\itfam\eightit}%
\def\csc{\eightcsc}%
\setbox\strutbox=\hbox{\vrule height7pt depth2pt width0pt}%
\normalbaselineskip=0,8\normalbaselineskip\normalbaselines\rm}
\newcount\absFnno\absFnno1
\write\AUX{\relFnno1}
\newif\ifMARKE\MARKEtrue
{\catcode`\@=11
\gdef\footnote{\ifMARKE\edef\@sf{\spacefactor\the\spacefactor}\/%
$^{\cite{Fn\the\absFnno}}$\@sf\fi
\MARKEtrue
\insert\footins\bgroup\eightpoint
\interlinepenalty100\let\par=\endgraf
\leftskip=0pt\rightskip=0pt
\splittopskip=10pt plus 1pt minus 1pt \floatingpenalty=20000\smallskip
\item{$^{\cite{Fn\the\absFnno}}$}%
\expandwrite\AUX{\neverexpand\ref{Fn\the\absFnno}{\neverexpand\the\relFnno}}%
\global\advance\absFnno1\write\AUX{\advance\relFnno1}%
\bgroup\strut\aftergroup\@foot\let\next}}
\skip\footins=12pt plus 2pt minus 4pt
\dimen\footins=30pc
\output={\plainoutput\immediate\write\AUX{\relFnno1}}
\newcount\Abschnitt\Abschnitt0
\def\beginsection#1. #2 \par{\advance\Abschnitt1%
\vskip0pt plus.10\vsize\penalty-250
\vskip0pt plus-.10\vsize\bigskip\vskip\parskip
\edef\TEST{\number\Abschnitt}
\expandafter\ifx\csname#1\endcsname\TEST\relax\else
\immediate\write16{#1 hat sich geaendert!}\fi
\expandwrite\AUX{\neverexpand\ref{#1}{\TEST}}
\leftline{\marginnote{#1}\bf\number\Abschnitt. \ignorespaces#2}%
\nobreak\smallskip\noindent\SATZ1\GNo0}
\def\Proof:{\par\noindent{\it Proof:}}
\def\Remark:{\ifdim\lastskip<\medskipamount\removelastskip\medskip\fi
\noindent{\bf Remark:}}
\def\Remarks:{\ifdim\lastskip<\medskipamount\removelastskip\medskip\fi
\noindent{\bf Remarks:}}
\def\Definition:{\ifdim\lastskip<\medskipamount\removelastskip\medskip\fi
\noindent{\bf Definition:}}
\def\Example:{\ifdim\lastskip<\medskipamount\removelastskip\medskip\fi
\noindent{\bf Example:}}
\def\Examples:{\ifdim\lastskip<\medskipamount\removelastskip\medskip\fi
\noindent{\bf Examples:}}
\newif\ifmarginalnotes\marginalnotesfalse
\newif\ifmarginalwarnings\marginalwarningstrue

\def\marginnote#1{\ifmarginalnotes\hbox to 0pt{\eightpoint\hss #1\ }\fi}

\def\strutdepth{\dp\strutbox}
\def\Randbem#1#2{\ifmarginalwarnings
{#1}\strut
\setbox0=\vtop{\eightpoint
\rightskip=0pt plus 6mm\hfuzz=3pt\hsize=16mm\noindent\leavevmode#2}%
\vadjust{\kern-\strutdepth
\vtop to \strutdepth{\kern-\ht0
\hbox to \hsize{\kern-16mm\kern-6pt\box0\kern6pt\hfill}\vss}}\fi}

\def\Zitat!{\Randbem{\bf?}{\bf Zitat}}

\newcount\SATZ\SATZ1
\def\proclaim #1. #2\par{\ifdim\lastskip<\medskipamount\removelastskip
\medskip\fi\goodbreak
\noindent{\bf#1.\ }{\it#2\Par}
\ifdim\lastskip<\medskipamount\removelastskip\goodbreak\medskip\fi}
\def\Aussage#1{\expandafter\def\csname#1\endcsname##1.{\resetitem
\ifx?##1?\relax\else
\edef\TEST{#1\penalty10000\ \number\Abschnitt.\number\SATZ}
\expandafter\ifx\csname##1\endcsname\TEST\relax\else
\immediate\write16{##1 hat sich geaendert!}\fi
\expandwrite\AUX{\neverexpand\ref{##1}{\TEST}}\fi
\proclaim {\marginnote{##1}\number\Abschnitt.\number\SATZ. #1\global\advance\SATZ1}.}}
\Aussage{Theorem}
\Aussage{Proposition}
\Aussage{Corollary}
\Aussage{Lemma}
\font\la=lasy10
\def\strich{\hbox{$\vcenter{\hbox
to 1pt{\leaders\hrule height -0,2pt depth 0,6pt\hfil}}$}}
\def\dashedrightarrow{\hbox{%
\hbox to 0,5cm{\leaders\hbox to 2pt{\hfil\strich\hfil}\hfil}%
\kern-2pt\hbox{\la\char\string"29}}}

\def\Bindestrich{\penalty10000-\hskip0pt}
\let\_=\Bindestrich
\def\.{{\sfcode`.=1000.}}

\def\Par{\par}
\def\:={\mathrel{\raise0,9pt\hbox{.}\kern-2,77779pt
\raise3pt\hbox{.}\kern-2,5pt=}}
\def\=:{\mathrel{=\kern-2,5pt\raise0,9pt\hbox{.}\kern-2,77779pt
\raise3pt\hbox{.}}} \def\mod{/\mskip-5mu/}
\def\into{\hookrightarrow}
\def\pfeil{\rightarrow}

\def\Pfeil{\longrightarrow}
\def\pf#1{\buildrel#1\over\rightarrow}
\def\Pf#1{\buildrel#1\over\longrightarrow}

\def\Ugleich{\hbox{$\cup$\kern.5pt\vrule depth -0.5pt}}
\def\|#1|{\mathop{\rm#1}\nolimits}
\def\<{\langle}
\def\>{\rangle}
\let\Times=\times
\def\times{\mathop{\Times}}
\let\Otimes=\otimes
\def\otimes{\mathop{\Otimes}}
\catcode`\@=11
\def\hex#1{\ifcase#1 0\or1\or2\or3\or4\or5\or6\or7\or8\or9\or A\or B\or
C\or D\or E\or F\else\message{Warnung: Setze hex#1=0}0\fi}
\def\fontdef#1:#2,#3,#4.{%
\alloc@8\fam\chardef\sixt@@n\FAM
\ifx!#2!\else\expandafter\font\csname text#1\endcsname=#2
\textfont\the\FAM=\csname text#1\endcsname\fi
\ifx!#3!\else\expandafter\font\csname script#1\endcsname=#3
\scriptfont\the\FAM=\csname script#1\endcsname\fi
\ifx!#4!\else\expandafter\font\csname scriptscript#1\endcsname=#4
\scriptscriptfont\the\FAM=\csname scriptscript#1\endcsname\fi
\expandafter\edef\csname #1\endcsname{\fam\the\FAM\csname text#1\endcsname}
\expandafter\edef\csname hex#1fam\endcsname{\hex\FAM}}
\catcode`\@=12 

\fontdef Ss:cmss10,,.
\fontdef Fr:eufm10,eufm7,eufm5.
\def\fa{{\Fr a}}

\def\fc{{\Fr c}}

\def\fg{{\Fr g}}
\def\fh{{\Fr h}}

\def\fl{{\Fr l}}

\def\fp{{\Fr p}}

\def\fs{{\Fr s}}
\def\ft{{\Fr t}}

\fontdef bbb:msbm10,msbm7,msbm5.
\fontdef mbf:cmmib10,cmmib7,.
\fontdef msa:msam10,msam7,msam5.
\def\CC{{\bbb C}}

\def\ZZ{{\bbb Z}}

\def\cU{{\cal U}}\def\cW{{\cal W}}
\def\cZ{{\cal Z}}
\mathchardef\leer=\string"0\hexbbbfam3F
\mathchardef\subsetneq=\string"3\hexbbbfam24
\mathchardef\semidir=\string"2\hexbbbfam6E
\mathchardef\dirsemi=\string"2\hexbbbfam6F
\mathchardef\haken=\string"2\hexmsafam78
\mathchardef\auf=\string"3\hexmsafam10
\def\Aq{{\overline{A}}}

\def\Cq{{\overline{C}}}


\def\mq{{\overline{m}}}
\def\nq{{\overline{n}}}

\def\sq{{\overline{s}}}\def\Sq{{\overline{S}}}

\def\Uq{{\overline{U}}}
\def\vq{{\overline{v}}}\def\Vq{{\overline{V}}}

\def\Xq{{\overline{X}}}

\def\zq{{\overline{z}}}
%
\newdimen\Parindent
\Parindent=\parindent

\def\textindent#1{\noindent\hskip\Parindent\llap{#1\enspace }\ignorespaces}
\def\itemitem{\par\indent\hangindent2\Parindent\textindent}


\abovedisplayskip 9.0pt plus 3.0pt minus 3.0pt
\belowdisplayskip 9.0pt plus 3.0pt minus 3.0pt
\newdimen\Grenze\Grenze2\Parindent\advance\Grenze1em
\newdimen\Breite
\newbox\DpBox
\def\NewDisplay#1
#2$${\Breite\hsize\advance\Breite-\hangindent
\setbox\DpBox=\hbox{\hskip2\Parindent$\displaystyle{%
\ifx0#1\relax\else\eqno{#1}\fi#2}$}%
\ifnum\predisplaysize<\Grenze\abovedisplayskip\abovedisplayshortskip
\belowdisplayskip\belowdisplayshortskip\fi
\global\futurelet\nexttok\WEITER}
\def\WEITER{\ifx\nexttok\qed\expandafter\leftQEDdisplay
\else\leftdisplay\fi}
\def\leftdisplay{\hskip-\hangindent\leftline{\box\DpBox}$$}
\def\leftQEDdisplay{\hskip-\hangindent
\line{\copy\DpBox\hfill\lower\dp\DpBox\copy\QEDbox}%
\belowdisplayskip0pt$$\bigskip\let\nexttok=}
\everydisplay{\NewDisplay}
\newcount\GNo\GNo=0
\newcount\maxEqNo\maxEqNo=0
\def\eqno#1{%
\global\advance\GNo1
\edef\FTEST{$(\number\Abschnitt.\number\GNo)$}
\ifx?#1?\relax\else
\ifnum#1>\maxEqNo\global\maxEqNo=#1\fi%
\expandafter\ifx\csname E#1\endcsname\FTEST\relax\else
\immediate\write16{E#1 hat sich geaendert!}\fi
\expandwrite\AUX{\neverexpand\ref{E#1}{\FTEST}}\fi
\llap{\hbox to 40pt{\ifx?#1?\relax\else\marginnote{E#1}\fi\FTEST\hfill}}}

\catcode`@=11
\def\eqalignno#1{\null\!\!\vcenter{\openup\jot\m@th\ialign{\eqno{##}\hfil
&\strut\hfil$\displaystyle{##}$&$\displaystyle{{}##}$\hfil\crcr#1\crcr}}\,}
\catcode`@=12

\newbox\QEDbox
\newbox\nichts\setbox\nichts=\vbox{}\wd\nichts=2mm\ht\nichts=2mm
\setbox\QEDbox=\hbox{\vrule\vbox{\hrule\copy\nichts\hrule}\vrule}
\def\qed{\leavevmode\unskip\hfil\null\nobreak\hfill\copy\QEDbox\medbreak}
\newdimen\HIindent
\newbox\HIbox
\def\setHI#1{\setbox\HIbox=\hbox{#1}\HIindent=\wd\HIbox}
\def\HI#1{\par\hangindent\HIindent\hangafter=0\noindent\leavevmode
\llap{\hbox to\HIindent{#1\hfil}}\ignorespaces}

\newdimen\maxSpalbr
\newdimen\altSpalbr
\newcount\Zaehler


\newif\ifxxx

{\catcode`/=\active

\gdef\beginrefs{%
\xxxfalse
\catcode`/=\active
\def/{\string/\ifxxx\hskip0pt\fi}
\def\TText##1{{\xxxtrue\tt##1}}
\expandafter\ifx\csname Spaltenbreite\endcsname\relax
\def\Spaltenbreite{1cm}\immediate\write16{Spaltenbreite undefiniert!}\fi
\expandafter\altSpalbr\Spaltenbreite
\maxSpalbr0pt
\gdef\alt{}
\def\\##1\relax{%
\gdef\neu{##1}\ifx\alt\neu\global\advance\Zaehler1\else
\xdef\alt{\neu}\global\Zaehler=1\fi\xdef\SigText{##1\the\Zaehler}}
\def\L|Abk:##1|Sig:##2|Au:##3|Tit:##4|Zs:##5|Bd:##6|S:##7|J:##8|xxx:##9||{%
\def\SigText{##2}\global\setbox0=\hbox{##2\relax}
\edef\TEST{[\SigText]}
\expandafter\ifx\csname##1\endcsname\TEST\relax\else
\immediate\write16{##1 hat sich geaendert!}\fi
\expandwrite\AUX{\neverexpand\ref{##1}{\TEST}}
\setHI{[\SigText]\ }
\ifnum\HIindent>\maxSpalbr\maxSpalbr\HIindent\fi
\ifnum\HIindent<\altSpalbr\HIindent\altSpalbr\fi
\HI{\marginnote{##1}[\SigText]}
\ifx-##3\relax\else{##3}: \fi
\ifx-##4\relax\else{##4}{\sfcode`.=3000.} \fi
\ifx-##5\relax\else{\it ##5\/} \fi
\ifx-##6\relax\else{\bf ##6} \fi
\ifx-##8\relax\else({##8})\fi
\ifx-##7\relax\else, {##7}\fi
\ifx-##9\relax\else, \TText{##9}\fi\Par}
\def\B|Abk:##1|Sig:##2|Au:##3|Tit:##4|Reihe:##5|Verlag:##6|Ort:##7|J:##8|xxx:##9||{%
\def\SigText{##2}\global\setbox0=\hbox{##2\relax}
\edef\TEST{[\SigText]}
\expandafter\ifx\csname##1\endcsname\TEST\relax\else
\immediate\write16{##1 hat sich geaendert!}\fi
\expandwrite\AUX{\neverexpand\ref{##1}{\TEST}}
\setHI{[\SigText]\ }
\ifnum\HIindent>\maxSpalbr\maxSpalbr\HIindent\fi
\ifnum\HIindent<\altSpalbr\HIindent\altSpalbr\fi
\HI{\marginnote{##1}[\SigText]}
\ifx-##3\relax\else{##3}: \fi
\ifx-##4\relax\else{##4}{\sfcode`.=3000.} \fi
\ifx-##5\relax\else{(##5)} \fi
\ifx-##7\relax\else{##7:} \fi
\ifx-##6\relax\else{##6}\fi
\ifx-##8\relax\else{ ##8}\fi
\ifx-##9\relax\else, \TText{##9}\fi\Par}
\def\Pr|Abk:##1|Sig:##2|Au:##3|Artikel:##4|Titel:##5|Hgr:##6|Reihe:{%
\def\SigText{##2}\global\setbox0=\hbox{##2\relax}
\edef\TEST{[\SigText]}
\expandafter\ifx\csname##1\endcsname\TEST\relax\else
\immediate\write16{##1 hat sich geaendert!}\fi
\expandwrite\AUX{\neverexpand\ref{##1}{\TEST}}
\setHI{[\SigText]\ }
\ifnum\HIindent>\maxSpalbr\maxSpalbr\HIindent\fi
\ifnum\HIindent<\altSpalbr\HIindent\altSpalbr\fi
\HI{\marginnote{##1}[\SigText]}
\ifx-##3\relax\else{##3}: \fi
\ifx-##4\relax\else{##4}{\sfcode`.=3000.} \fi
\ifx-##5\relax\else{In: \it ##5}. \fi
\ifx-##6\relax\else{(##6)} \fi\PrII}
\def\PrII##1|Bd:##2|Verlag:##3|Ort:##4|S:##5|J:##6|xxx:##7||{%
\ifx-##1\relax\else{##1} \fi
\ifx-##2\relax\else{\bf ##2}, \fi
\ifx-##4\relax\else{##4:} \fi
\ifx-##3\relax\else{##3} \fi
\ifx-##6\relax\else{##6}\fi
\ifx-##5\relax\else{, ##5}\fi
\ifx-##7\relax\else, \TText{##7}\fi\Par}
\bgroup
\baselineskip12pt
\parskip2.5pt plus 1pt
\hyphenation{Hei-del-berg Sprin-ger}
\sfcode`.=1000
\beginsection References. References

}}

\def\endrefs{%
\expandwrite\AUX{\neverexpand\ref{Spaltenbreite}{\the\maxSpalbr}}
\ifnum\maxSpalbr=\altSpalbr\relax\else
\immediate\write16{Spaltenbreite hat sich geaendert!}\fi
\egroup\write16{Letzte Gleichung: E\the\maxEqNo}
\write16{Letzte Aufzaehlung: I\the\maxItemcount}}



\fontdef sans:cmss10,,.

\def\sp{{\fs\fp}}
\def\reg{_{\rm reg}}

\let\UL=\underline
\newdimen\uldepth

\def\underline#1{\setbox0=\hbox{$#1$}\uldepth=\dp0%
\setbox0=\hbox{$\UL{#1}$}\dp0=\uldepth\box0}

\let\nOL=\overline
\newdimen\olheight

\def\overline#1{\setbox0=\hbox{$#1$}\olheight=\ht0%
\setbox0=\hbox{$\nOL{#1}$}\ht0=\olheight\box0}

\def\half{{\textstyle{1\over2}}}
\def\fgq{{\overline\fg}}
\def\thetaq{{\overline\theta}}
\def\sigmaq{{\overline\sigma}}
\def\faq{{\overline\fa}}
\def\flq{{\overline\fl}}

\input xy
\xyrequire{curve}
\xyrequire{cmtip}
\xyrequire{matrix}
\xyrequire{arrow}
\SelectTips{cm}{10}

\newdir{ (}{{}*!/-5pt/\dir{^(}}
\def\inj{\ar@{^(->}}

\def\cxymatrix#1{\vcenter{\xymatrix@=15pt{#1}}}


\title{Invariant functions on symplectic representations}

\beginsection Introduction. Introduction

Let $G$ be a connected reductive group over $\CC$. Our goal is to
study invariants of a {\it symplectic representation}, i.e., a finite
dimensional $G$\_representations $V$ which is equipped with a
non\_degenerate symplectic form $\omega$. Our main tool will be the
$\fg^*$\_valued covariant
$$
m:V\pfeil\fg^*:v\mapsto\ell_v\quad
{\rm where}\quad \ell_v(\xi):=\half\omega(\xi v,v)
$$
called the {\it moment map} (here, $\fg^*$ is the coadjoint
representation of $G$). The key idea is to construct (very special)
invariants on $V$ by pulling back the (well\_known) invariants on $\fg^*$.

This construction can be described more geometrically as follows. It
is known (Chevalley) that $G$\_invariants on $\fg^*\cong\fg$ are in
bijection with $W_G$\_invariants on $\ft^*$ (where $\ft^*$ is the dual
of a Cartan subspace and $W_G$ is the Weyl group). Thus, we can define
the composed morphism
$$
m\mod G:V\Pf m\fg^*\Pfeil\ft^*/W_G
$$
called the {\it invariant moment map}. Then every function on
$\ft^*/W_G$ gives rise to an invariant function on $V$.

Let $R_0\subseteq\CC[V]^G$ be the ring of invariants obtained this
way.  It turns out to be more natural to consider the (slightly)
larger ring $R$ of invariants which are algebraic over $R_0$. Then our
main results are:

\item{$\bullet$}$R$ is a polynomial ring. More precisely, there is a
subspace $\fa^*\subseteq\ft^*$ and a reflection group $W_V$ acting on
$\fa^*$ such that
$$
R=\CC[\fa^*]^{W_V}.
$$
The dimension of $\fa^*$ is called the {\it symplectic rank} of
$V$. The group $W_V$ is a subtle invariant of $V$ called its {\it
little Weyl group}.

\item{$\bullet$}Both the ring of all functions $\CC[V]$ and the ring
of all invariants $\CC[V]^G$ are free $R$\_modules. Geometrically,
this means: let
$$31
\xymatrix@R=15pt@C=0pt{
V\ar[dr]_(0.3)\psi\ar[rr]^\pi&&V\mod G\ar[dl]^(0.3){\psi\mod G}
\\
&\kern-10pt\fa^*/W_V\kern-10pt&
}
$$
be the commutative triangle induces by the inclusions
$\CC[\fa^*]^{W_V}=R\subseteq k[V]^G\subseteq \CC[V]$.  Then both $\psi$ and
$\psi\mod G$ are faithfully flat. In particular, all fibers have the
same dimension.

\item{$\bullet$}The fibers of $\psi\mod G$ are called the {\it
symplectic reductions} of $V$. They form a flat family of
$2c$\_dimensional Poisson varieties. Here $c$ is the {\it symplectic
complexity} of $V$. We show that a generic symplectic reduction contains
a dense open subset which is isomorphic to an open subset of
$\CC^{2c}$ with its standard symplectic structure.

\item{$\bullet$} The generic fibers of $\pi:V\pfeil V\mod G$ are fiber
products of the form $G\times^LF$ where $L$ is a Levi subgroup and
$$
F=A\times\CC^{2m_1}\times\ldots\times\CC^{2m_s}.
$$
Here $A$ is a torus and $L$ acts on $F$ via a surjective homomorphism
$$
L\auf A\times Sp_{2m_1}(\CC)\times\ldots\times Sp_{2m_s}(\CC).
$$

\noindent Of particular interest is the case when $R=\CC[V]^G$, i.e.,
when the moment map furnishes ``almost all'' invariants. These
representations are called {\it multiplicity free}. We prove:

\item{$\bullet$}The following are equivalent (see
\S\cite{multiplicityfree} for unexplained terminology):

\itemitem{---} $V$ is multiplicity free.

\itemitem{---} All invariants Poisson\_commute with each other.

\itemitem{---} The algebra of invariants $\cW(V)^G$ inside the Weyl
algebra $\cW(V)$ is commutative.

\itemitem{---} The generic $G$\_orbits are coisotropic.\Par

\item{$\bullet$} Every multiplicity free symplectic representation is
cofree, i.e., $\CC[V]$ is a free module over $\CC[V]^G$ (this
condition implies that $\CC[V]^G$ is a polynomial ring).

\noindent The cofreeness property of multiplicity free representations
is very restrictive. We used it in \cite{mfclass} to classify all
multiplicity free symplectic representations.

If $U$ is any finite dimensional representation of $G$ then $V=U\oplus
U^*$ carries a canonical symplectic structure. In fact, $V$ can be
considered as the cotangent bundle of $U$. Therefore, the present
paper can be seen as an extension of my theory of invariants on
cotangent bundles, started in \cite{WuM}. This also explains some
terminology: $V=U\oplus U^*$ is multiplicity free in the symplectic
sense if and only if $\CC[U]$ is multiplicity free in the usual sense.

The main tool for studying the geometry of a cotangent bundle was the
local structure theorem of Brion\_Luna\_Vust \cite{BLV}. In this paper
we prove a symplectic analog of the structure theorem. More
precisely, we construct a Levi subgroup $M\subseteq G$ and an
$M$\_stable subspace $S\subseteq V$ such that
\item{$\bullet$}The restriction of $\omega$ to $S$ is
non\_degenerate.

\item{$\bullet$}There is a dense open subset $S^0\subseteq S$ and an
embedding $q:S^0\into V$ such that the diagram
$$32
\cxymatrix{
S^0\ar[r]^q\ar[d]_{m_S\mod M}&V\ar[d]^{m_V\mod G}
\\
\ft^*/W_M\ar[r]&\ft^*/W_G
}
$$
commutes. One key point is that $q$ is not the natural inclusion of
$S^0$ in $V$. In general, it will not be even linear. Another key
point is that $S^0$ is not just any open subset: it is
the complement of an explicitly given hyperplane. This will imply that
$S^0$ meets the zero\_fiber of $m_S\mod M$ . This is vital for proving
the equidimensionality of $m_V\mod G$ by induction on $\|dim|V$.

\item{$\bullet$} The map
$$
G\times^MS^0\pfeil V:[g,s]\mapsto gq(s)
$$
is \'etale on an open  dense subset (we are actually proving something
much more precise). This statement links the generic structure of $V$
with the generic structure of $S$. It implies, in particular, that
$S^0$ and $V$ have essentially the same image in $\ft^*/W_G$ (see
diagram \cite{E32}).

\item{$\bullet$}The morphism $S^0\mod M\pfeil V\mod G$ induced by $q$
is a Poisson morphism. This will be used to study the generic
symplectic reductions.

\medskip

\noindent{\bf Acknowledgment:} Most of this paper was written during a
year\_long stay at the University of Freiburg, Germany. I would like
to thank this institution for its hospitality.

\medskip

\noindent {\bf Notation:} 1. We are working in the category of complex
algebraic varieties even though $\CC$ could be replaced by any
algebraically closed field of characteristic zero. The ring of regular
functions on a variety $X$ is denoted by $\CC[X]$ while $\CC(X)$ is
its field of rational functions. The dual of a vector space $V$ will
be denoted by $V^*$. In contrast, the one\_dimensional torus is
$\CC^\Times=\CC\setminus\{0\}$.

2. If $H$ is any algebraic group and $X$ an affine $H$\_variety we
denote the categorical quotient by $X\mod H$, i.e., $X\mod
H=\|Spec|\CC[X]^H$.

3. The Lie algebra of any group is denoted by the corresponding
fraktur letter. If $X$ is an $H$\_variety then every $\xi\in\fh$
induces a vector field $\xi_*$ on $X$. Its value at $x\in X$ is
denoted by $\xi x$.

4. In the whole paper, $G$ will denote a connected reductive group. We
sometimes use tacitly that $\fg^*\cong\fg$ as $G$\_varieties.  We
choose a Borel subgroup $B\subseteq G$ and a maximal torus $T\subseteq
B$.  The root system of $\fg$ is denoted by $\Delta$, the positive
roots by $\Delta^+$. In any root subspace $\fg_\alpha$ we choose a
generator $\xi_\alpha$ such that
$[\xi_\alpha,\xi_{-\alpha}]=\alpha^\vee$.

\beginsection momentmap. The moment map

An {\it affine Poisson variety} is an affine variety $X$ equipped with a
bilinear map
$$
\CC[X]\times\CC[X]\Pfeil\CC[X]:(f,g)\mapsto\{f,g\},
$$
called the {\it Poisson product}, satisfying:

\item{1.} the Poisson product is a Lie algebra structure on $\CC[X]$
and
\item{2.} for each $f\in\CC[X]$, the map $g\mapsto\{f,g\}$ is a
derivation of $\CC[X]$.

\noindent Because of the second condition, the concept of a Poisson
variety is local in nature. In particular, every open subset of $X$
(in the Zariski- or in the \'etale sense) is again a Poisson variety.

A morphism $\pi:X\pfeil Y$ between two Poisson varieties is
a {\it Poisson morphism} if
$$
\{f,g\}_Y=\{f\circ\pi,g\circ\pi\}_X\quad f,g\in\CC[Y].
$$

There are two main examples of Poisson varieties:

\noindent1. {\it Dual Lie algebras:} Let $\fg$ be a finite dimensional
Lie algebra and $X=\fg^*$. Since $\CC[X]=S^*\fg$ we have
$\fg\subseteq\CC[X]$. One can show that the Lie bracket on $\fg$ extends
uniquely to a Poisson product on $\CC[X]$.

\noindent2. {\it Symplectic varieties:} Let $X$ be a smooth variety
$X$ equipped with a non\_degenerate closed $2$-form $\omega$. Then the
Poisson structure is constructed as follows. First, every function
$f\in\CC[X]$ gives rise to a $1$\_form $df$. Since $\omega$ is
non\_degenerate, it identifies the tangent bundle of $X$ with its
cotangent bundle. Thus, $df$ corresponds to a vector field $H_f$,
called the {\it Hamiltonian vector field} attached to $f$. The
relationship of $f$ and $H_f$ is expressed in the formula
$$
\omega(\xi,H_f)=\xi(f)=df(\xi)
$$
where $\xi$ is any tangent vector. Then
$$
\{f,g\}:=\omega(H_f,H_g)=-df(H_g)=-H_g(f)=dg(H_f)=H_f(g).
$$
defines a Poisson product. Observe the formula
$H_{\{f,g\}}=[H_f,H_g]$. Therefore, the map
$$
\CC[X]\Pfeil\Gamma(T_X):f\mapsto H_f
$$
is a Lie algebra homomorphism (with kernel $\CC$).

Now assume that the algebraic group $G$ acts on the symplectic variety
$X$ in such a way that $\omega$ is $G$\_stable. Then every $\xi\in\fg$
induces a vector field $\xi_*$ on $X$. Moreover, the map
$\xi\mapsto-\xi_*$ is a Lie algebra homomorphism\footnote{The minus
sign comes from the fact that $\xi$ acts on functions naturally from
the right while a left action is used to define the bracket of vector
fields.}. Thus we obtain the diagram
$$22
\cxymatrix{&\fg\ar[d]^{-\xi_*}\ar@{.>}[ld]_{m^*}\\
\CC[X]\ar[r]_{H_f}&\Gamma(T_X)}
$$
The symplectic $G$\_variety is called {\it Hamiltonian} if it is
equipped with a Lie algebra homomorphism $m^*:\fg\pfeil\CC[X]$ such
that diagram \cite{E22} commutes. More
geometrically, we consider the morphism
$$
m:X\pfeil\fg^*: x\mapsto[\xi\mapsto m^*(\xi)(x)]
$$
called the {\it moment map} of $X$. It is also characterized by the
equations
$$0
\eqalignno{
21&\omega(\xi x,\xi' x)&=\<m(x),[\xi,\xi']\>
\quad\hbox{for all }\xi,\xi'\in\fg\cr
18&\omega(\xi x,\eta)&=\<D_xm(\eta),\xi\>
\quad\hbox{for all }\xi\in\fg,\eta\in T_xX.\cr}
$$
Here, \cite{E21} expresses the fact that $m^*$ is a Lie algebra
homomorphism while \cite{E18} says that \cite{E22} commutes. These two
equations also imply that $m$ is $\fg$\_equivariant. Since $G$ is
assumed to be connected we conclude that $m$ is even
$G$\_equivariant. Given the $G$\_action and the symplectic structure
of $X$, the moment map is unique up to translation by an element of
$(\fg^*)^G$.  In fact, if $m'$ is another moment map then \cite{E18}
implies that $m-m'$ has everywhere a zero derivative. Observe the
immediate consequence of \cite{E18}:
$$19
\|ker|D_xm=(\fg x)^\perp\quad{\rm and}\quad\|Im|D_xm=\fg_x^\perp
$$

Equation \cite{E21} implies that $m$ is a Poisson morphism.
Moreover, a function $f\in\CC[X]$ is $G$\_invariant if and only if
$$
\{m^*(\xi),f\}=H_{m^*(\xi)}(f)=-\xi_*(f)=0
$$
for all $\xi\in\fg$. Thus, $\CC[X]^G$ is the Poisson centralizer of
$m^*\CC[\fg^*]$ in $\CC[X]$. This implies, in particular, that
$m^*\CC[\fg^*]^G$ is in the Poisson center of $\CC[X]^G$. The
geometric counterpart to this is the {\it invariant moment map}
$$
m\mod G:X\pfeil\fg^*\mod G\cong\ft^*/W_G
$$
which is the composition of the moment map $X\pfeil\fg^*$ with the
categorical quotient map $\fg^*\pfeil\ft^*/W_G$. Not only is this a
Poisson map (with the trivial Poisson bracket on $\ft^*/W_G$) but also
all fibers are Poisson varieties.

In this paper we are mostly concerned with the special case that $X=V$
is a finite dimensional $G$\_representation equipped with a
non\_degenerate 2-form $\omega\in\wedge^2V^*$. Then there is a
canonical moment map namely
$$
m:V\pfeil\fg^*:v\mapsto[\xi\mapsto\half\omega(\xi v,v)].
$$
This implies the following useful formula: let $\fh\subseteq\fg$ be a
subalgebra and $\fh^\perp\subseteq\fg^*$ its annihilator. Then
$$
m^{-1}(\fh^\perp)=
\{v\in V\mid\omega(\xi v,v)=0\ {\rm for\ all}\ \xi\in\fh\}.
$$

\beginsection local. The symplectic local structure theorem

In this section we develop our main technical tool, a symplectic
version of the local structure theorem of Brion\_Luna\_Vust,
\cite{BLV}.

Let $V$ be a symplectic $G$\_representation. The structure theorem
will depend on the choice of a highest weight vector $v_0\in V$. The
construction will not work for the defining representation of the
symplectic group.

\Definition: Let $U\subset V$ be the submodule $U$ generated by
$v_0$. Then $U$ (or $v_0$) is called {\it singular} if $U$ is an
anisotropic subspace of $V$ (hence itself symplectic) and $G\pfeil
Sp(U)$ is surjective. A dominant weight $\chi$ of $G$ is {\it
singular} if the corresponding irreducible representation $U$ is
singular.

\Remarks: 1. If the root system of $G$ has a component of type ${\Ss
C}_n$ (for any $n\ge1$, including ${\Ss C}_1={\Ss A}_1$) then the
first fundamental weight of that component is singular. Conversely,
all singular dominant weights are of this form.

2. Clearly, if $U\subseteq V$ is singular then its highest weight is
singular. The converse is false since the singularity of $U$ depends
on its embedding in $V$. Let, for example, $G=Sp_{2n}(\CC)$ and
$V=\CC^{2n}\oplus\CC^{2n}$. Then both summands are singular, but the
submodule $\{(v,iv\mid v\in\CC^n\}$ is isotropic, hence non\_singular.

\Lemma 2chi. Let $U$ be an irreducible symplectic $G$\_module with
highest weight $\chi$. Then $\chi\not\in\Delta^+$. Moreover, if one
of the following conditions holds then $\chi$ is singular:
\Item{1} $2\chi\in\Delta^+$.
\Item{6} $2\chi=2\alpha-\beta$ with $\alpha,\beta\in\Delta^+$.
\Item{2} $2\chi=\alpha+\beta$ with $\alpha,\beta\in\Delta^+$.

\Proof: We first prove that each of the conditions
\cite{I1}--\cite{I2} implies that $\chi$ is singular.

\cite{I1} Assume $2\chi=\alpha\in\Delta^+$. Then $\half\alpha$ is a
weight and $\<\alpha,\beta^\vee\>$ is even for all
$\beta\in\Delta^+$. The same holds for the simple root in the
$W$\_orbit of $\alpha$. Thus, $\alpha$ is root in a root system of
type ${\Ss C}_n$, $n\ge1$ and $\chi={1\over2}\alpha$ is its first
fundamental weight.

\cite{I6} Since $\half\beta$ is a weight, we are as above in a root
system of type ${\Ss C}_n$, $n\ge1$. There, the assertion can be
checked directly.

\cite{I2} We claim that if $\chi=\half(\alpha+\beta)$ is a dominant
weight then either $2\chi$ is a root (and we are done by\cite{I1}) or the
corresponding simple module is not symplectic. First observe, that
$$
\<\chi,\alpha^\vee\>=1+\half\<\beta,\alpha^\vee\>\in\ZZ,\quad
\<\chi,\beta^\vee\>=1+\half\<\alpha,\beta^\vee\>\in\ZZ
$$
implies that both $\<\alpha,\beta^\vee\>$ and $\<\beta,\alpha^\vee\>$
are even integers. This is only possible if both of them are
zero. This implies
$$23
\<\chi,\alpha^\vee\>=1,\quad \<\chi-\alpha,\beta^\vee\>=1.
$$
Let $L\subseteq G$ be the smallest Levi subgroup having $\alpha$ and
$\beta$ as roots. Its semisimple rank is at most two. Let
$U'\subseteq U$ be the $L$\_module generated by a highest weight
vector of $U$. Then \cite{E23} implies that $U'$ contains also the
lowest weight vector of $U$. Therefore, the restriction of the
invariant symplectic form of $U$ to $U'$ is non\_zero. We conclude
that $U'$ is an irreducible symplectic $L$\_module. Thus we are
reduced to $\|rk|G=2$ where the claim is easily
checked case\_by\_case.

Finally, if $\chi$ were a root since then \cite{I2} implies that
$\chi$ were singular. Thus $\chi\not\in\Delta^+$.\qed

{\it For the rest of this section we assume that $v_0$ is
non\_singular}. Let $P$ be the stabilizer of the line $\CC
v_0$. Denote its unipotent radical by $P_u$ and its Levi subgroup by
$M$. Let $P^-=MP_u^-$ be the opposite parabolic subgroup. Let $v_0^-\in
V$ be a lowest weight vector with $\omega(v_0^-,v_0)=1$. Its weight is
necessarily $-\chi$. We are going to need the following observation:

\Lemma. The subspace $\CC v_0\oplus\fp_u^-v_0\subseteq V$ is
isotropic.

\Proof: Suppose there are root vectors
$\xi_{-\alpha},\xi_{-\beta}\in\fp_u^-$ with
$\omega(\xi_{-\alpha}v_0,\xi_{-\beta}v_0)\ne0$. Then
$(\chi-\alpha)+(\chi-\beta)=0$ which implies that $\chi$ is a singular
weight (\cite{2chi}). Since $v_0$ is non\_singular, even $\<Gv_0\>$ would be
isotropic. The same argument shows $\omega(\fp_u^-v_0,v_0)=0$.\qed

Now consider the following subspace of $V$:
$$
S:=(\fp_u^-v_0)^\perp\cap(\fp_uv_0^-)^\perp=\{v\in V\mid
\omega(\fp_u^-v_0,v)=0,\omega(\fp_uv_0^-,v)=0\}.
$$
It is obviously a representation under $M$. We will see later (proof
of \cite{S0S0}) that the restriction of $\omega$ to $S$ is
non\_degenerate. Thus $S$ it is a symplectic $M$\_representation. The goal
of this section is to describe a reduction procedure from $(V,G)$ to
$(S,M)$.

The natural inclusion of $S$ into $V$ is not compatible with invariant
moment maps. To achieve compatibility, we alter it in a non\_linear
fashion. More precisely, we look at the subset
$$
\Sigma:=(\fp_u^-v_0)^\perp\cap m^{-1}(\fp_u^\perp)=\{v\in V\mid
\omega(\fp_u^-v_0,v)=0,\omega(\fp_uv,v)=0\}
$$
whose definition is very similar to that of $S$. 

\Definition: For any subset $Z\subseteq V$ we put $Z^0:=\{v\in Z\mid
\omega(v,v_0)\ne0\}$.

\medskip

\noindent Now we have:

\Lemma projection. There is a decomposition
$V=\fp_u^-v_0\oplus(\fp_uv_0^-)^\perp$ and the projection to the
second summand induces an isomorphism $p:\Sigma^0\pf\sim S^0$.

\Proof: We claim that $\omega$ induces a perfect pairing between
$\fp_u^-v_0$ and $\fp_uv_0^-$. Indeed,
$\omega(\xi_{-\alpha}v_0,\xi_\beta v_0^-)=0$ for $\alpha\ne\beta$
while $\omega(\xi_{-\alpha}v_0,\xi_\alpha
v_0^-)=-\omega([\xi_\alpha,\xi_{-\alpha}]v_0,v_0^-)=\<\chi|\alpha^\vee\>\ne0$.
The claim implies $\fp_u^-v_0\cap(\fp_uv_0^-)^\perp=0$, hence
$V=\fp_u^-v_0\oplus(\fp_uv_0^-)^\perp$ for dimension reasons.

Now $\fp_u^-v_0\subseteq(\fp_u^-v_0)^\perp$ implies that $p$ maps
$(\fp_u^-v_0)^\perp$ onto $S$. In particular, $p(\Sigma)\subseteq
S$. Moreover $\omega(\fp_u^-v_0,v_0)=0$ implies $p(\Sigma^0)\subseteq
S^0$. We are going to construct the inverse map. More precisely, we
construct a morphism $\phi:S^0\pfeil\fp_u^-$ such that the inverse map
is given by $s\mapsto s+\phi(s)v_0$.

Let $s\in S^0$, $\xi_-\in\fp_u^-$ and put $v=s+\xi_-v_0$. Then
$v\in\Sigma$ means $\omega(\xi_+v,v)=0$ for all $\xi_+\in\fp_u$.
We have
$$1
0=\omega(\xi_+v,v)=\omega(\xi_+s,s)+\omega(\xi_+s,\xi_-v_0)+
\omega(\xi_+\xi_-v_0,s)+\omega(\xi_+\xi_-v_0,\xi_-v_0).
$$
The last summand vanishes by \cite{2chi}\cite{I6}. Because of
$\omega(\xi_+s,\xi_-v_0)=\omega(\xi_+\xi_-v_0,s)$, equation \cite{E1}
becomes
$$16
\omega(\xi_+\xi_-v_0,s)=-\half\omega(\xi_+s,s)
\quad\hbox{for all $\xi_+\in\fp_u$.}
$$
Clearly, it suffices to show that for any $s\in S^0$ these equations
have a unique solution.  First observe that \cite{E16} is a square
system of inhomogeneous linear equations for $\xi_-$. Thus, it
suffices to show that the matrix
$$
\big[\omega(\xi_\alpha\xi_{-\beta}v_0,s)\big]_{\alpha,\beta\in\Delta_u}
$$
is invertible (where $\Delta_u$ is the set of roots $\alpha$ with
$\xi_\alpha\in\fp_u$). We claim that it is even triangular with
non\_zero diagonal entries. Let $U:=\<Gv_0\>$ and
$U^-:=\<Gv_0^-\>$. These are two (not necessarily distinct) simple
$G$\_modules and one is the dual of the other. We have $U^-\cap
U^\perp=0$, and therefore $V=U^-\oplus U^\perp$. Since $v_0$ and
$\xi_\alpha\xi_{-\beta} v_0$ are elements of $U$, we may replace $s$ by
its component in $U^-$. Then $s$ has a weight decomposition
$s=\sum_\eta s_\eta$ with $\eta\ge-\chi$ (meaning $\eta+\chi$ being a
sum of positive roots). Assume
$\omega(\xi_\alpha\xi_{-\beta}v_0,s)\ne0$. Then
$\eta=-\chi+\beta-\alpha\ge-\chi$, hence $\beta\ge\alpha$. This shows
that the matrix is triangular. For the diagonal terms with
$\beta=\alpha$ we get
$$
\omega(\xi_\alpha\xi_{-\alpha}v_0,s)=
\omega([\xi_\alpha,\xi_{-\alpha}]v_0,s)=
\<\chi|\alpha^\vee\>\omega(v_0,s)\ne0.
$$\qed

The sought-after embedding of $S^0$ into $V^0$ is the composition
$$
\cxymatrix{q:S^0\ar[r]_{p^{-1}}^\sim&\ \Sigma^0\ \inj[r]&V^0}
$$

\Lemma S0S0. Both subsets $S$ and $\Sigma^0$ inherit the structure of
a Hamiltonian $M$\_variety from $V$. Moreover, the map
$q:S^0\Pf\sim\Sigma^0$ is an isomorphism of Hamiltonian
$M$\_varieties. In particular, the following diagram commutes:
$$35
\cxymatrix{S^0\ar[r]^q\ar[d]_{m_S}&V\ar[d]^{m_V}
\\
\fl^*&\fg^*.\ar[l]^{\rm res}
}
$$

\Proof: First, we have $S\cap
S^\perp\subseteq(\fp_uv_0^-)^\perp\cap\fp_u^-v_0=0$
(\cite{projection}). This shows that $S$ is anisotropic and therefore
a symplectic subspace of $V$. The set $\Sigma^0$ is smooth since it
is, via $p$, isomorphic to an open subset of $S$.

Fix $s_1\in S^0$ and put $v_1:=q(s_1)\in\Sigma^0$. Then we get a
map of tangent spaces $Dq:T_{s_1}S^0\pfeil T_{v_1}\Sigma^0$. Let
$s,\sq\in T_{s_1}S^0=S$ be tangent vectors and consider their
images $v=Dq(s)$, $\vq=Dq(\sq)$ in $T_{v_1}\Sigma^0$. Now recall that
$$
q(s)=s+\phi(s)v_0
$$
where $\phi:S^0\pfeil\fp_u^-$ is some morphism. Thus, there are
$\xi_-,\overline\xi_-\in\fp_u^-$ such that $v=s+\xi_-v_0$,
$\vq=\sq+\overline\xi_-v_0$. From this we get
$$
\omega(v,\vq)=\omega(s,\sq)+\omega(\xi_-v_0,\sq)+
\omega(s,\overline\xi_-v_0)+\omega(\xi_-v_0,\overline\xi_-v_0)
=\omega(s,\sq)
$$
since $S\subseteq(\fp_u^-v_0)^\perp$ and $\fp_u^-v_0$ is
isotropic. This shows that $\omega|_{\Sigma^0}$ is non\_degenerate and
that $q$ is a symplectomorphism.

We also have $v_1=s_1+\xi_1^-v_0$ for some $\xi_1^-\in\fp_u^-$. Let
$\xi\in\fl$. Then
$$
\eqalign{2 m_{\Sigma^0}(v_1)(\xi)&=\omega(\xi v_1,v_1)=\cr
&=\omega(\xi s_1,s_1)+
\omega(\xi s_1,\xi_1^-v_0)+
\omega(\xi\xi_1^-v_0 ,s_1)+
\omega(\xi\xi_1^-v_0,\xi_1^-v_0)\cr
&=\omega(\xi s_1,s_1)=2 m_S(s_1)(\xi).}
$$
since $\xi s_1\in S$ and $\xi\xi_1^-v_0\in\fp_u^-v_0$. This shows the
commutativity of \cite{E35}.\qed

Note that the bottom arrow in \cite{E35} goes in the ``wrong''
direction. This is fixed by using the invariant moment maps:

\Theorem Commute. The following diagram commutes:
$$15
\cxymatrix{
S^0\ar[r]^q\ar[d]_{m_S\mod M}&V\ar[d]^{m_V\mod G}
\\
\ft^*/W_M\ar[r]&\ft^*/W_G
}
$$

\Proof: We use $\fg^*\cong\fg$. Then $\fl^*\cong\fl$,
$\ft^*\cong\ft$, and $\fg^*\pfeil\fl^*$ becomes the orthogonal
projection $\fg\pfeil\fl$. The definition of $\Sigma$ implies
$(m_V\circ q)(S^0)=m_V(\Sigma^0)\subseteq\fp_u^\perp=\fp$. In view of
\cite{S0S0}, the assertion now follows from the commutativity of
$$
\cxymatrix{
\fl\ar[d]&\fp\ar@{>>}[l]\ar[d]
\\
\ft/W_M\ar[r]&\ft/W_G.
}
$$
One can see that this diagram commutes observe, e.g., that
$\fp\pfeil\fl$ is the categorical quotient by $P_u$ (see
e.g. \cite{E4} below).\qed

The embedding $q$ is, in general, not a Poisson morphism. We just have
a statement for $P_u$\_invariants:

\Theorem. The map $q\mod P_u:S^0\pfeil V^0\mod P_u$ is a Poisson morphism.

\Proof: In view of \cite{S0S0}, we may replace $S^0$ by
$\Sigma^0$. For $f\in\CC[V^0]$ let $f_\Sigma$ be its restriction to
$\Sigma^0$. For all $f,g\in\CC[V^0]^{P_u}$ we have to show that
$\{f,g\}=\{f_\Sigma,g_\Sigma\}$.

Fix a point $s\in\Sigma^0$. Let $H\in T_sV^0$ be the tangent
vector which is dual to $dg$. It is defined by
$$11
\omega(H,\eta)=dg(\eta)\quad\hbox{for all $\eta\in T_sV^0$.}
$$
Then we have $\{f,g\}(s)=df(H)$. Similarly, let $H_\Sigma\in
T_s\Sigma^0$ be dual to $dg_\Sigma$, i.e., with
$$12
\omega(H_\Sigma,\eta)=dg(\eta)\quad\hbox{for all $\eta\in T_s\Sigma^0$.}
$$
The composition $V\pf m\fg^*\pfeil\fp_u^*$ is the moment map for
$P_u$. Its zero fiber is $Z:=m^{-1}(\fp_u^\perp)$, hence $s\in
Z$. From \cite{E19} we get $T_sZ=(\fp_us)^\perp$. In particular, $Z$
is smooth in $s$ if its isotropy group inside $P_u$ is trivial. Thus,
we get from \cite{PSigma} that
$$
(\fp_us)^\perp=T_sZ=\fp_us\oplus T_s\Sigma^0.
$$
From $H_\Sigma\in T_s\Sigma^0\subseteq(\fp_us)^\perp$ we infer
$\omega(H_\Sigma,\eta)=0$ for all $\eta\in\fp_us$. On the other hand, the
$P_u$\_invariance of $g$ implies $dg(\fp_us)=0$. Therefore,
\cite{E12} is valid for all $\eta\in\fp_us\oplus
T_s\Sigma^0=(\fp_us)^\perp$. Combined with \cite{E11} we get
$H_\Sigma-H\in(\fp_su)^{\perp\perp}=\fp_us$. Hence,
$$
\{f_\Sigma,g_\Sigma\}(s)=df(H_\Sigma)=df(H)+df(H_\Sigma-H)=\{f,g\}(s).
$$
For the last equality, we used the $P_u$\_invariance of $f$.\qed

Restricting to $G$\_invariants yields:

\Corollary poissonL. The map $q\mod G:S^0\mod M\pfeil V\mod G$ is a
Poisson morphism.

Finally, we connect the generic geometry of $S$ with that of $V$. We do
that in two steps.

\resetitem

\Lemma PSigma. There is an isomorphism
$$
P\times^M\Sigma^0\pf\sim m^{-1}(\fp_u^\perp)^0.
$$

\Proof: Consider the subspace $N=(\fp_u v_0^-)^\perp\subseteq
V$. Then the local structure theorem of Brion\_Luna\_Vust, \cite{BLV}
Prop.~1.2, implies that
$$
P\times^MN^0\pfeil V^0
$$
is an isomorphism. Intersecting both sides with the $P$\_stable subset
$ m^{-1}(\fp_u^\perp)=\{v\in V\mid\omega(\fp_u v,v)=0\}$ yields the
assertion.\qed

Let $W_G$ and $W_M$ be the Weyl groups of $G$ and $M$,
respectively. The morphism $\ft^*/W_M\pfeil\ft^*/W_G$ is smooth in an
open subset denoted by $(\ft^*/W_M)\reg$. For every morphism
$Z\pfeil\ft^*/W_M$ let $Z\reg$ be the preimage of $(\ft^*/W_M)\reg$.
For example, $\ft^*\reg$ is the set of points $\chi\in\ft^*$ with
$\<\chi|\alpha^\vee\>\ne0$ for all roots corresponding to $\fp_u$. In
particular, it is an affine variety. Its most important property is:
$$3
\xi\in\ft\reg\Longrightarrow C_G(\xi)\subseteq M\ {\rm and}\
C_\fg(\xi)\subseteq\fl.
$$
If $Z$ is a Hamiltonian $M$\_variety with moment map $Z\pfeil\fl^*$
then we can form $Z\reg$ with respect to the map
$Z\pfeil\fl^*\pfeil\fl^*\mod M=\ft^*/W_M$.

We continue with two well\_known lemmas:

\Lemma AuxLemma0. The map
$$4
P\times^M\fl\reg\Pf\sim\fp\reg.
$$
is an isomorphism.

\Proof: We start with a general remark concerning fiber products. Let
$G$ be a group, $H\subseteq G$ a subgroup, $Y$ an $H$\_variety, $X$ a
$G$\_variety, and $Y\subseteq X$ an $H$\_invariant
subvariety. Consider the induced map $\Phi:G\times^HY\pfeil X$. Then
$\Phi$ is surjective if and only if

\item{$(S)$} $GY=X$.

\noindent Moreover, $\Phi$ is injective if and only if

\item{$(I)$} $g\in G$, $y\in Y$, $gy\in Y$ imply $g\in H$.

We use this criterion to show that \cite{E4} is bijective. Let
$\xi\in\fp\reg$ with semisimple part $\xi_s$. After conjugation with
$P$ we may assume $\xi_s\in\fl\reg$. But then $\xi\in
C_\fg(\xi_s)\subseteq\fl$. This shows $(S)$. Now let $u\in P_u$,
$\xi\in\fl\reg$ with $u\xi\in\fl\reg$. Since $u\xi=\xi$ modulo $\fp_u$
we have $u\xi=\xi$. Thus $u\in C_P(\xi_s)\subseteq M$, hence
$u=1$. This shows $(I)$. Thus, \cite{E4} is a bijective morphism
between normal varieties and therefore an isomorphism (\cite{LuAO}
Lemme~1.8).\qed

\Lemma AuxLemma. The following diagram is Cartesian:
$$2
\cxymatrix{
G\times^P\fp\reg\ar[r]\ar[d]&\fg\ar[d]\\
(\ft/W_M)\reg\ar[r]&\ft/W_G
}
$$

\Proof: Because of \cite{E4} it suffices to show that the diagram
$$5
\cxymatrix{
G\times^M\fl\reg\ar[r]\ar[d]&\fg\ar[d]\\
(\ft/W_M)\reg\ar[r]&\ft/W_G
}
$$
is Cartesian or, in other words, that
$$30
G\times^M\fl\reg\pfeil\tilde\fg:=\fg\times\limits_{\ft/W_G}(\ft/W_M)\reg
$$
is an isomorphism. We follow the same strategy as for the proof of
\cite{AuxLemma0}. 

For (S) let $(\xi,W_M\xi')\in\tilde\fg$, i.e., $\xi\in\fg$, $\xi'\in\ft\reg$,
and $\xi,\xi'$ have the same image in $\ft/W_G$. Then there is $g\in
G$ such that $g\xi_s=\xi'$. Hence $g\xi\in C_\fg(\xi')=\fl$. Thus,
$(\xi,W_M\xi')$ is the image of $[g^{-1},g\xi]$.

For (I) assume $g\in G$ and $\xi\in\fl\reg$ are elements with
$g\xi\in\fl\reg$ and $\xi$, $g\xi$ have the same image in
$\ft/W_M$. Thus, there is $l\in M$ with $lg\xi_s=\xi_s$. Since
$C_G(\xi_s)=M$ we conclude $g\in M$.

Finally, $\tilde\fg$ is an \'etale cover of
$\fg$, hence smooth. Thus, Luna's Lemma implies that \cite{E30} is an
isomorphism.\qed

\noindent The next statement provides the link between the geometry of
$S$ and that of $V$.

\Theorem LSThm. Let $X:=m_V^{-1}(\fp_u^\perp)\reg$. Then there is a
commutative diagram
$$25
\cxymatrix{
G\times^MS^0\reg\inj[r]^{\iota}\ar[dr]_{m_S\mod M}&
G\times^PX\ar[r]\ar[d]&V\ar[d]^{m_V\mod G}
\\
&(\ft^*/W_M)\reg\ar[r]&\ft^*/W_G
}
$$
where the right hand square is Cartesian and where $\iota$ is an open
embedding. Moreover, $S^0\reg\ne\emptyset$.

\Proof: We use the identifications $\fg^*=\fg$, $\ft^*=\ft$, and
$\fp_u^\perp=\fp$. Now consider the diagram
$$
\cxymatrix{
G\times^P m^{-1}(\fp\reg)\ar[r]\ar[d]&V\ar[d]^m
\\
G\times^P\fp\reg\ar[r]\ar[d]&\fg\ar[d]
\\
(\ft/W_M)\reg\ar[r]&\ft/W_G
}
$$
The top square is clearly Cartesian while the bottom square is
Cartesian by \cite{AuxLemma}. This shows that the right square of
\cite{E25} is Cartesian. The commutativity of the left square follows
from \cite{Commute}. The morphism $\iota$ is induced by $q$ and is an
open embedding by \cite{PSigma}. Finally, consider $v=v_0+v_0^-\in S$. Then
$$
 m_S(v)(\xi)=\half\omega(\xi v,v)=
\cases{0&if $\xi\in\fp_u+\fp_u^-$\cr-\chi(\xi)&if $\xi\in\fl$.\cr}
$$
This shows $ m_S(v)\in S^0\reg\ne\emptyset$.\qed

\beginsection TermReps. Terminal representations

Clearly the reduction step of the preceding section doesn't work if
all highest weight vectors are singular. But also a highest weight
vector which spans a 1\_dimensional submodule has to be avoided since
in that case $M=G$ and $S^0=V^0$ and the local structure theorem
becomes tautological.

\Definition: A highest weight vector of $V$ is called {\it terminal}
if it is either singular or generates a one\_dimensional $G$\_module
(i.e., $M=G$). The representation $V$ is {\it terminal} if all of its
highest weight vectors are terminal. A dominant weight is called {\it
terminal} if it is singular or a character of $\fg$.

\medskip\noindent Here is the classification:

\Proposition terminal. Let $V$ be a terminal symplectic
representation. Then $(G,V)$ decomposes as
$$
\eqalign{V=&V_0\oplus V_1\oplus\ldots\oplus V_s\cr
G=&G_0\times Sp(V_1)\times\ldots\times Sp(V_s)\cr}
$$
where $V_0$ is a direct sum of 1\_dimensional $G_0$\_modules.

\Proof: Let $V=V_0\oplus\Vq$ where $V_0$ is the sum of all
1\_dimensional subrepresentations. Let $V_1\subseteq\Vq$ be a simple
submodule and let $V'$ be its complement in $\Vq$. Then $V_1$ is
anisotropic and $G\pfeil Sp(V_1)$ is surjective. Since the symplectic
group is simply connected, there is a unique splitting
$G=Sp(V_1)\times G'$. We claim that $Sp(V_1)$ acts trivially on
$V'$. Indeed, otherwise it would contain a simple submodule
$V_1'$. Since also that submodule is singular, there is an isomorphism
$\phi:V_1\pf\sim V_1$. But then the submodule $\{v+i\phi(v)\mid v\in
V_1\}$ of $V_1\oplus V_1'$ would be isotropic, hence
non\_singular. This shows that $(G,\Vq)=(Sp(V_1),V_1)\times(G',V')$
and we are done by induction.\qed

Now we analyze the moment map of terminal representations. First the
non\_toric part:

\Lemma symplecticcase. Let $V=\CC^{2n}$ be a symplectic vector space and
$G=Sp_{2n}(\CC)$. Then, using the identification
$\sp_{2n}(\CC)^*\cong\sp_{2n}(\CC)$, the moment map maps $V$ to the set of
nilpotent matrices of rank $\le1$. In particular, the invariant moment
map $V\pfeil\ft^*/W_G$ is zero.

\Proof: Assume the symplectic structure is given by the skewsymmetric
matrix $J$, i.e., $\omega(u,v)=u^tJv$. Let $A\in\sp_{2n}(\CC)$. Then
$$
2m(v)(A)=\omega(Av,v)=-v^tJAv=-\|tr|(v^tJAv)=\|tr|((-vv^tJ)A).
$$
Thus $m(v)=-\half vv^tJ$ which has rank one. Moreover
$\|tr|vv^tJ=\omega(v,v)=0$. Hence $m(v)$ is nilpotent.\qed

Next, we deal with the toric part:

\Lemma toruscase. Let $A$ be a torus acting with a finite kernel on
a symplectic representation $V$. Let $m:V\pfeil\fa^*$ be the
corresponding moment map. Let $C$ be an irreducible component of
$m^{-1}(0)$. Then there is a section $\sigma:\fa^*\pfeil V$ with
$\sigma(\fa^*)\cap C\ne\emptyset$.

\Proof: There is a direct sum decomposition
$$17
V=\oplus_{i=1}^n(\CC_{\chi_i}\oplus\CC_{-\chi_i})
$$ where $\CC_\chi$ denotes the one\_dimensional representation of
$\fa=\|Lie|A$ for the character $\chi\in\fa^*$. Let $x_i,y_i$
($i=1,\ldots,n$) be the corresponding coordinates. The moment map for
$A$ acting on $V$ is $m_V(x_i,y_i)=\sum_{i=1}^nx_iy_i\chi_i$. We will
call $\chi_i$ {\it critical} if $\chi_i$ is not in the span of
$\{\chi_j\mid j\ne i\}$.

We prove the assertion by induction on $\|dim|V$. First, assume that
none of the characters are critical. Then we claim that the zero fiber
of $m_V$ is irreducible. To see this, we factor $m_V$ into
$$
\mq:V\pfeil\CC^n:(x_i,y_i)\mapsto(x_iy_i)\quad{\rm
and}\quad\pi:\CC^n\pfeil\fa^*:(t_i)\mapsto\sum_{i=1}^nt_i\chi_i.
$$
The map $\pi$ is linear. Let $K$ be its kernel.  Then
$m_V^{-1}(0)=\mq^{-1}(K)$. Since $\mq$ is visibly faithfully flat, the
same holds for $\mq^{-1}(K)\pfeil K$. Therefore, every irreducible
component maps dominantly to $K$, i.e., we just have to check that the
generic fiber is irreducible. The non\_existence of a critical weight
means precisely that every coordinate of a generic point of $K$ is
non\_zero. Therefore, the fiber is isomorphic to $(\CC^\Times)^n$,
hence irreducible. This proves the claim.

Since the action is locally faithful, the characters $\chi_i$ span
$\fa^*$, Thus, we may arrange that $\{\chi_1,\ldots,\chi_m\}$
forms a basis of $\fa^*$. Then the set $\fa^0$ defined
by the equations
$$
x_1=\ldots=x_m=1,x_{m+1}=\ldots=x_n=y_1=\ldots=y_n=0
$$
is subset of $V$ which maps isomorphically onto $\fa^*$. In
particular, it hits the zero fiber $C$.

Now assume that $\chi_1$, say, is critical. Let $\fa':=\CC\chi_1$,
$V':=\CC_{\chi_1}\oplus\CC_{-\chi_1}$ and
$\fa'':=\<\chi_2,\ldots,\chi_n\>$,
$V'':=\oplus_{i=2}^n\CC_{\chi_i}\oplus\CC_{-\chi_i}$. Then we have
decompositions $V=V'\oplus V''$ and
$\fa^*=\fa'\oplus\fa''$, and it suffices to prove the assertion for
$(V',\fa')$ and $(V'',\fa'')$ separately. For the latter case we use
the induction hypothesis. In the former case, $m^{-1}(0)$ has two
components which intersect either the section $\{1\}\times\CC$ or
$\CC\times\{1\}$.\qed

Now we use the structure of terminal representations and the reduction
procedure of \S\cite{local} to give a preliminary description of the
generic structure of $V$. This will be refined further below
(\cite{generic3}).

\Lemma generic2. For every symplectic $G$\_representation $V$ there is a
commutative diagram of $G$\_varieties
$$27
\cxymatrix{
(G\times^LF)\times U\ar[r]\ar[d]_p&V\ar[d]^{m_V\mod G}
\\
\ft^*/W_L\ar[r]&\ft^*/W_G
}
$$
which identifies $(G\times^LF)\times U$ with an open subset of the
fiber product $V\times_{\ft^*/W_G}\ft^*/W_L$. Here:
\item{$\bullet$}$L$ is a Levi subgroup of $G$ equipped with a
surjective homomorphism
$$
\phi:L\auf A\times Sp_{2m_1}(\CC)\times\ldots\times Sp_{2m_s}(\CC)
$$
where $A$ is a torus and $s$ can be $0$.
\item{$\bullet$}$F=A\times \CC^{2m_1}\times\ldots\times\CC^{2m_s}$
with $L$\_action induced by $\phi$.
\item{$\bullet$}$U$ is a non\_empty open subset of $\fa^*\times\CC^{2c}$.
\item{$\bullet$}$p$ is the composition $(G\times^LF)\times U\auf
U\pfeil\fa^*\pfeil\ft^*/W_L$.\Par
\noindent
Further properties are $U\reg\ne\leer$ and the induced morphism on
$G$\_invariants $U\pfeil V\mod G$ is a Poisson morphism (where
$\CC^{2s}$ has the standard symplectic structure).


\Proof: We construct $L$, $F$, and $U$ by induction on $\|dim|V$.
Assume first, that $V$ is terminal. Then it has a decomposition
$$
V=V_0\oplus V_1\oplus\ldots\oplus V_s
$$
as in \cite{terminal}. By \cite{symplecticcase}, the invariant moment map
of $V$ factors through $V_0$. Since $V_0$ is a symplectic
representation, there is a $G$\_stable decomposition
$V_0=\oplus_{i=1}^n(\CC_{\chi_i}\oplus\CC_{-\chi_i})$ as in
\cite{E17}. Then $V_0=U_0\oplus U_0^*$ with
$U_0=\oplus_{i=1}^n\CC_{\chi_i}$. Consider the torus
$C:=(\CC^\Times)^n$. Then the map
$$26
C\times\fc^*=(\CC^\Times)^n\times(\CC^n)^*\Pfeil U_0\oplus U_0^*:
(u_i,v_i)\mapsto (u_i,u_i^{-1}v_i)
$$
realizes the cotangent bundle $T^*C$ with its natural symplectic
structure as a dense open subset of $U_0\oplus U_0^*$. The group $G$
acts on $V_0$ via a homomorphism $G\pfeil C$. Let $A\subseteq C$ be
its image. Since $A$ is connected it has a complement $\Aq$ in $C$,
i.e., we get a splitting $C=A\times\Aq$. Using \cite{E26}, we get a
$G$\_equivariant open embedding
$$
T^*A\times T^*\Aq=T^*C\into V_0
$$
which is compatible with the symplectic structure and where $G$ acts
trivially on $\Aq$. Using \cite{E26} for $\Aq$ instead of $C$ we see
that $T^*\Aq$ is isomorphic to an open subset of $\CC^{2c}$,
$c=\|dim|\Aq$, with its standard symplectic structure. Thus, we have
proved the existence of the diagram \cite{E27} with $L=G$, $F=A\times
V_1\times\ldots\times V_s$, and $U=\fa^*\times T^*\Aq$.

If $V$ is not terminal then from \cite{LSThm} we get a diagram
$$34
\cxymatrix{
G\times^MS^1\ar[r]\ar[d]&V\ar[d]
\\
(\ft^*/W_M)\reg\ar[r]&\ft^*/W_G
}
$$
(with $S^1:=S^0\reg$) which identifies the upper left corner with an
open subset of the fiber product. Since $V$ is not terminal we may
assume $S\ne V$. Thus, the induction hypothesis furnishes a diagram
$$33
\cxymatrix{
(M\times^LF)\times U\ar[r]^(0.7)\alpha\ar[d]&S\ar[d]
\\
\ft^*/W_L\ar[r]&\ft^*/W_M
}
$$
Since $S_1$ is an {\it affine} open subset of $S$, its preimage under
$\alpha$ is so as well. This implies that we may shrink $U$ such that,
in diagram{E33},we may replace $S$ by $S^1$. Now replace $\alpha$ by
its fiber product $G\times^M\alpha$ and compose the ensuing diagram
with diagram \cite{E34} yielding diagram \cite{E27} with required
properties.

By induction, we know that the image $\Uq$ of $U$ in $\ft^*/W_M$
intersects the locus over which $\ft^*/W_L\pfeil\ft^*/W_M$ is
unramified. From the fact that $S^1\ne\leer$ we get that $\Uq$
contains points in which $\ft^*/W_M\pfeil\ft^*/W_G$ is
unramified. Since $\Uq$ is irreducible it contains points where both
conditions hold showing $U\reg\ne\leer$.  Finally, the statement about
the Poisson structure follows from \cite{poissonL}.  \qed

In the proof of \cite{section1} and \cite{haesslich} we need to deal
with one\_dimensional submodules directly. For that we use the
following reduction:

\Lemma onedimensional. Assume $v_0\in V$ generates a one\_dimensional
submodule with character $\chi$. Let $v_0^-\in V$ be a
$G$\_eigenvector with $\omega(v_0^-,v_0)=1$. Put $\Vq:=(\CC
v_0\oplus\CC v_0^-)^\perp$ and $\fgq:=\|ker|\chi\subseteq\fg$. Let
$\mq:\Vq\pfeil\fgq^*$ be the moment map for $\fgq$. Then there is a
$G$\_equivariant isomorphism $\Phi$ such that the following diagram
commutes:
$$
\cxymatrix{
\Vq\times\CC\times\CC^\Times\ar[d]_{m_0}\ar[r]^(0.68)\Phi&V^0\ar[d]^m
\\
\fgq^*\oplus\CC\chi\ar[r]^(0.62)\sim&\fg^*
}
$$
Here $m_0$ is the map $(v,t,y)\mapsto(\mq(v),ty)$. Moreover, $\fg$ acts
trivially on $\CC$ and with the character $-\chi$ on $\CC^\Times$.

\Proof: If $\chi=0$ then we define $\Phi(v,t,y)=v+tv_0+yv_0^-$. So
assume $\chi\ne0$ from now on. Choose $\xi$ in the center of $\fg$
with $\chi(\xi)=1$. This yields splittings $\fg=\fgq\oplus\CC\xi$
and $\fg^*=\fgq^*\oplus\CC\chi$. We define a $G$\_invariant
function on $\Vq$ by
$$
f(v)=\<m(v),\xi\>=\half\omega(\xi v,v).
$$
Then $m(v+xv_0+yv_0^-)=\mq(v)+(f(v)+xy)\chi$ and
$$
\Phi(v,t,y):=v+{t-f(v)\over y}v_0+yv_0^-
$$
has the required property.\qed

\ignore

Finally, we construct $U_0$. Put $X:=G\times^LF$. Then $X$ contains a
dense $G$\_orbit. Thus, taking $G$\_invariants, we get from $X\times
U\pfeil V$ a morphism $U\pfeil V\mod G$.  Let $\tilde
V:=V\times_{\ft^*/W_G}\ft^*/W_L$. Then $\tilde U:=(G\times^LF)\times
U$ is an open subset of $\tilde V$. By shrinking $U$, we may assume
that $\tilde U$ is affine. Then 

\Corollary closure. The closure of the image of $V$ in $\ft^*/W_G$ is
equal to the image of $\fa^*$ in $\ft^*/W_G$.

\medskip

The moment map $m:V\pfeil\fg^*$ induces a map $m\mod G:V\mod
G\pfeil\fg^*\mod G$ on categorical quotients. The non\_empty fibers of
this morphism are called the {\it symplectic reductions} of
$V$. The next statement is well\_known in the context of compact group
actions (see e.g. \cite{GS}).

\Theorem weakreduction. Let $V$ be a symplectic representation. Then
every symplectic reduction carries a Poisson structure. Moreover, for
a generic symplectic reduction there is an open dense subset on which
the Poisson structure is symplectic, i.e., non\_degenerate.
;;

\beginsection equi. Equidimensionality

Next, we use the results of section~\cite{local} to investigate the
geometry of the invariant moment map. In abuse of notation ($W_G$ does not, in
general, act on $\fa^*$) we denote the image of $\fa^*$ in $\ft^*/W_G$
by $\fa^*/W_G$.  

\Theorem section1. Let $V$ be a symplectic representation and let
$\fa^*\subseteq\ft^*$ be as in \cite{generic2}. Then:
\Item{30} The image of $V$ in $\ft^*/W_G$ is $\fa^*/W_G$.
\Item{31} There is a morphism $\sigma:\fa^*\pfeil V$ such that
$$14
\cxymatrix{
\fa^*\ar[r]^\sigma\ar@{>>}[rd]&V\ar[d]^{m\mod G}
\\
&\fa^*/W_G
}
$$
is commutative.
\Item{32} Let $C$ be an irreducible component of $(m\mod G)^{-1}(0)$. Then
one can choose $\sigma$ such that $\sigma(\fa^*)\cap C\ne\emptyset$.

Before we prove the theorem, we mention the following
consequence:

\Corollary equidim. The invariant moment map $V\pfeil\fa^*/W_G$ is
equidimensional.

\Proof: By general properties of fiber dimensions, every fiber has
dimension $\ge\|dim|V-\|dim|\fa^*$. We have to show the opposite
inequality. By homogeneity and semicontinuity, the fiber of maximal
dimension is the zero fiber. Let $C$ be an irreducible component and
let $\sigma$ be as in \cite{section1}\cite{I31} and \cite{I32}. Then
$$
0=\|dim|\sigma(\fa^*)\cap C\ge\|dim|\fa^*+\|dim|C-\|dim|V
$$
shows indeed $\|dim|C\le\|dim|V-\|dim|\fa^*$.\qed

\noindent{\it Proof of \cite{section1}}: \cite{generic2} implies that
the image of $V$ in $\ft^*/W_G$ is contained in
$\fa^*/W_G$. Therefore, part \cite{I30} will follows from
part~\cite{I31}. To show \cite{I31} and \cite{I32}
we use induction on $\|dim|V$.

If $V$ is terminal then the assertions are covered by
\cite{symplecticcase} and \cite{toruscase}. Otherwise, the component
$C$ contains a non\_terminal highest weight vector $v_1$. We defer the
proof to \cite{haesslich} since it is quite technical.  Then
$v_0^-:=w_0v_1$ (with $w_0\in W_G$ the longest element) is a lowest
weight vector contained in $C$ and there is a highest weight vector
$v_0$ (possibly not in $C$) with $\omega(v_0^-,v_0)=1$.  Now we use
the notation of section~\cite{local} and consider the commutative
diagram~\cite{E15}. Since $v_0^-\in C$, the preimage $C_0$ of $C$ in
$S^0$ is not empty.  Then we have
$$
\|dim|\fa^*\ge\|codim|_VC\ge\|codim|_SC_0\ge
\|codim|_Sm_S^{-1}(0)\ge\|dim|\fa^*
$$
The last inequality holds because, by induction, $S\pfeil\fa^*/W_M$ is
equidimensional. Thus, we have equality throughout. This implies in
particular that every component of $C_0$ is a component of
$m_S^{-1}(0)$. Let $C_1$ be one of them. From \cite{onedimensional}
applied to $(S,M,v_0)$ we obtain a diagram
$$
\cxymatrix{
\Sq\times\CC\times\CC^\Times\ar[r]^(0.68)\sim\ar[d]_{m_0}&S^0\ar[d]^m
\\
\flq^*\oplus\CC\chi\ar[r]^(0.62)\sim&\fl^*
}
$$
In particular, $C_1\cong\Cq\times\CC^\Times$ where $\Cq$ is a component of
$\mq^{-1}(0)\subseteq\Sq$. By induction, there is a section
$\sigmaq:\faq^*\pfeil\Sq$ meeting $\Cq$. Thus we get a section
$\fa^*=\faq^*\oplus\CC\chi\pfeil S^0$ defined by
$$
\alpha+t\chi\mapsto\Phi(\sigmaq(\alpha),t,1)
$$
which meets $C_1$ and therefore $C_0$. Finally, the composition
$\sigma$ with the inclusion $S^0\pfeil V$ has the required
properties.\qed

The following lemma was used in the preceding proof.

\Lemma haesslich. Let $C$ be an irreducible component of
$(m\mod G)^{-1}(0)$ and assume $V$ is not terminal. Then $C$ contains a
non\_terminal highest weight vector.

\Proof: Let $U\subseteq V$ be the submodule generated by $C$. Assume
first, that $U$ has a non\_terminal highest weight. We claim that
there is a 1\_parameter subgroup $\rho:\G_m\pfeil T\subseteq G$ with
the following properties:
$$6
\<\alpha|\rho\>>0\quad\hbox{for all positive roots
$\alpha$}.
$$
$$7
\<\chi_1|\rho\><\<\chi_2|\rho\>\quad\hbox{for all $\chi_1$ terminal,
$\chi_2$ non\_terminal highest weights of $V$.}
$$
$$8
\<\chi_1|\rho\>\ne\<\chi_2|\rho\>\quad\hbox{for all weights $\chi_1\ne\chi_2$
of $V$.}
$$
An explicit calculation for the group $Sp_{2n}$ shows that \cite{E6}
and \cite{E7} hold for $\rho$ equal to the sum of the fundamental
coweights. Now a slight perturbation of $\rho$ yields additionally
\cite{E8}.

Let $\Gamma$ be the set of weights of $U$. Choose $v\in C$ such that
its $\chi$\_component $v_\chi$ is non\_zero for every $\chi\in\Gamma$.
By the last condition \cite{E8} on $\rho$, the maximum $N$ of
$\{\<\chi|\rho\>\mid\chi\in\Gamma\}$ is attained for a unique weight
$\chi_0$. Because of \cite{E6}, this weight is a highest
weight. Moreover, since $U$ contains a non\_terminal highest weight,
the weight $\chi_0$ is non\_terminal (condition \cite{E7}). Thus, the
limit
$$
v_0:=\|lim|_{t\pfeil\infty}t^{-N}\rho(t)v
$$
exists, is a non\_terminal highest weight vector (condition
\cite{E8}), and is contained in $C$ since $C$, as a component of
$(m\mod G)^{-1}(0)$, is a cone.

Now we are reduced to the case where all highest weights of $U$ are
terminal. If $U$ contains a one\_dimensional submodule then
\cite{onedimensional} reduces the statement to $\Vq$ and we are done
by induction. Thus $U$ contains only irreducible submodules $U_1$ such
that the image of $G$ in $GL(U_1)$ is the full symplectic
group. Assume that there is such a $U_1$ and it is isotropic. Then
there is a second isotropic submodule $U_2$ which is $G$\_isomorphic
to $U_1$ and such that $\Uq:=U_1\oplus U_2$ is anisotropic. The moment
map $m$ is invariant under the $\G_m$\_action which is
$(tu_1,t^{-1}u_2)$ on $\Uq$ and trivial on $\Uq^\perp$. Hence also $C$
is $\G_m$\_stable. Now choose $v\in C$ which has a non\_zero component
in $U_1$. Then
$$
v_1:=\|lim|_{t\pfeil\infty}t^{-1}(t\cdot v)
$$
exists, is non\_zero, and lies in $C\cap U_1$. Thus there is $g\in G$
such that $gv_1$ is a non\_singular highest weight vector in $C$.

Now we are reduced to the case that $U=U_1\oplus\ldots\oplus U_s$
where each $U_i$ is anisotropic and such that $G\pfeil
Sp(U_1)\times\ldots\times Sp(U_s)$ is surjective. Let
$Q:=U^\perp$. Then $V=U\oplus Q$. Clearly, we may assume that the
action of $G$ on $V$ is (locally) effective. Assume $Sp(U_i)\subseteq G$ acts
trivially on $Q$. Then the moment map $m_V$ factors through $V/U_i$
and we are done by induction. Thus, we may assume that even the action
of $G$ on $Q$ is locally effective. From $G\subseteq Sp(Q)$ we infer
$\|rk|G\le\half\|dim|Q$. On the other hand, $C$ has codimension
$\le\|dim|\ft^*=\|rk|G$ in $V$. From $C\subseteq U$ we get 
$$
\|dim|Q=\|codim|_VU\le\|codim|_VC\le
\|rk|G\le\half\|dim|Q.
$$
This implies that $V=U$ is terminal.\qed

\beginsection littleweylgroup. The little Weyl group

Next, we analyze the equidimensional morphism $V\pfeil\fa^*/W_G$. Let
$N(\fa^*)$ and $C(\fa^*)$ be the normalizer and the centralizer of
$\fa^*$ in $W_G$, respectively. Then the quotient $N(\fa^*)/C(\fa^*)$
acts faithfully on $\fa^*$.

\Theorem construction. Let $\fa^*\subseteq\ft^*$ be as in
\cite{generic2}. Then there is subgroup $W_V$ of $N(\fa^*)/C(\fa^*)$
and a morphism $\psi:V\pfeil\fa^*/W_V$ such that the generic fibers of
$\psi$ are irreducible and such that the following diagram commutes:
$$20
\cxymatrix{
V\ar[d]_\psi\ar[dr]^{m\mod G}
\\
\fa^*/W_V\ar[r]&\fa^*/W_G
}
$$
Moreover, the pair $(\fa^*,W_V)$ is, up to conjugation by $W_G$, uniquely
determined by this property.

\Proof: First we claim that, up to isomorphism, there is a unique
diagram
$$28
\cxymatrix{
V\ar[d]_\psi\ar[rd]^{m\mod G}
\\
Z\ar[r]^(0.4)\beta&\fa^*/W_G
}
$$
such that $Z$ is connected and normal, $\beta$ is finite, and the
generic fibers of $\psi$ are connected.  The morphism
$V\pfeil\fa^*/W_G$ is surjective. Moreover, since $\beta$ is finite,
$\psi$ has at least to be dominant. Thus we get inclusions
$$
\CC[\fa^*/W_G]\subseteq \CC[Z]\subseteq\CC[V].
$$
The requirement that the generic fibers of $\psi$ are connected means
that $\CC(Z)$ is algebraically closed in $\CC(V)$. This shows that
$\CC(Z)$ is the algebraic closure of $\CC(\fa^*/W_G)$ in
$\CC(V)$. Since $Z$ is normal and finite over $\fa^*/W_G$ we see that
$\CC[Z]$ is the integral closure of $\CC[\fa^*/W_G]$ in
$\CC(Z)$. Combined we get that $\CC[Z]$ is the integral closure of
$\CC[\fa^*/W_G]$ in $\CC(V)$ which proves uniqueness. For existence,
we take $Z=\|Spec|R$ where $R$ is the $\CC(V)$. This integral closure
is contained in $\CC[V]$ since $V$ is normal. Thus we get the required
factorization \cite{E28}.

It rests to identify $Z$ with $\fa^*/W_V$. For that consider the
morphism $\sigma:\fa^*\pfeil V$ from \cite{section1}\cite{I31}. Then
we have a commutative diagram
$$
\cxymatrix{
\fa^*\ar[d]\ar[dr]
\\
Z\ar[r]&\fa^*/W_G
}
$$
and therefore
$\CC[\fa^*/W_G]\subseteq\CC[Z]\subseteq\CC[\fa^*]$. Let
$\Gamma:=N(\fa^*)/C(\fa^*)$. Then the morphism $\fa^*\auf\fa^*/W_G$
factors through $\fa^*/\Gamma$ and we claim that
$$24
\fa^*/\Gamma\auf\fa^*/W_G
$$
is the normalization morphism. In fact, since it is clearly finite and
surjective and since $\fa^*/\Gamma$ is normal, we have to show that
\cite{E24} is birational. Suppose $w\in W_G$ with $w\not\in
N(\fa^*)$. Then $w\fa^*\ne\fa^*$ and therefore $w\xi\not\in\fa^*$ for
$\xi$ in an open dense subset of $\fa^*$. Since $W_G$ is finite, there
is a dense open subset $U$ of $\fa^*$ such that
$W_G\xi\cap\fa^*=\Gamma\xi$ for all $\xi\in U$. Thus \cite{E24} is
injective on $U$, proving the claim.

From the claim and the fact that $Z$ is normal, we get
$$
\CC[\fa^*]^\Gamma\subseteq\CC[Z]\subseteq\CC[\fa^*].
$$
For the field of fraction we obtain 
$$
\CC(\fa^*)^\Gamma\subseteq\CC(Z)\subseteq\CC(\fa^*).
$$
From Galois theory we get a unique subgroup $W_V$ of
$\Gamma$ with $\CC(Z)=\CC(\fa^*)^{W_V}$. The normality of $Z$ implies
$$
\CC[Z]=\CC(Z)\cap\CC[\fa^*]=\CC(\fa^*)^{W_V}\cap\CC[\fa^*]=
\CC[\fa^*]^{W_V}=\CC[\fa^*/W_V].
$$
This implies $Z=\fa^*/W_V$.

For the uniqueness statement observe that the preimage of $m\mod G(V)$
in $\ft^*$ is the union of the subspaces $w\fa^*$, $w\in W_G$. Thus,
$\fa^*$ is unique up to conjugation by $W_G$. Now assume also
$\fa^*/W_V'$ fits into the diagram \cite{E20}. Since $Z$ is unique,
there would be a $\fa^*/\Gamma$\_isomorphism
$\nq:\fa^*/W_V\pfeil\fa^*/W_V'$. Since $\fa^*$ is Galois over
$\fa^*/\Gamma$, the isomorphism $\nq$ extends to an automorphism
$n\in\Gamma$ of $\fa^*$. This shows $W_V'=nW_Vn^{-1}$.\qed

\Remark: The subspace $\fa^*$ determines the Levi subgroup $L$ of
\cite{generic2}. In fact, it follows from $U\reg\ne\leer$ in
\cite{generic2} that $L$ is the pointwise centralizer of $\fa^*$ (and
not larger than that). Thus, the roots of $L$ are those $\alpha$ with
$\<\fa^*,\alpha^\vee\>\ne0$. On the other hand, there is no canonical
parabolic $P$ with Levi subgroup $L$. More precisely, after repeated
application of the construction in section~\cite{local} one winds up
with a parabolic subgroup with Levi part $L$ but that subgroup would
depend on the choice of highest weights. A most trivial example is
$G=GL(n)$ with $n\ge3$ and $V=\CC^n\oplus(\CC^n)^*$. This
representation has two highest weight vectors and the corresponding
parabolic subgroups are not conjugated. Nevertheless, their Levi parts
are.

\Theorem mainresult. The subgroup $W_V$ of $GL(\fa^*)$ is generated by
reflections. In particular, $\fa^*/W_V$ is smooth.

\Proof: It follows from \cite{equidim} and the finiteness of
$\fa^*/W_V\pfeil\fa^*/W_G$ that also $\psi:V\pfeil\fa^*/W_V$ is
equidimensional. From the proof of \cite{construction} we know that
$\CC[\fa^*/W_V]$ is integrally closed in $\CC[V]$. Thus, the assertion
follows from the following lemma, essentially due to
Panyushev~\cite{Pan}.\qed

\Lemma. Let $X$ and $Y$ be vector spaces, $\Gamma\subseteq GL(Y)$ a
finite subgroup and $\psi:X\pfeil Y/\Gamma$ a dominant equidimensional
morphism. Assume that $\CC[Y/\Gamma]$ is integrally closed in
$\CC[X]$. Then $\Gamma$ is a generated by reflections and $Y/\Gamma$
is smooth.

\Proof: Let $\Gamma_r\subseteq\Gamma$ be the normal subgroup generated
by reflections. Then, by the Shepherd\_Todd\_Chevalley theorem
(\cite{Bou}\ Ch.V \S5 ${\rm n}^{\underline{\rm o}}$~5.5) the quotient
$Y/\Gamma_r$ is again an affine space. Moreover, $\Gamma/\Gamma_r$
acts linearly on it and does not contain any reflections. Thus we may
replace $Y$ and $\Gamma$ by $Y/\Gamma_r$ and $\Gamma/\Gamma_r$,
respectively. Thereby we achieve that $\Gamma$ does not contain any
reflections anymore and we have to show $\Gamma=1$.

Let $\Xq:=X\times_{Y/\Gamma}Y$. Since $Y\pfeil Y/\Gamma$ is
surjective, also $\Xq\pfeil X$ is surjective. Thus, we can find an
irreducible component $\tilde X$ of $\Xq$ such that $\tilde X\pfeil X$
is dominant. Thus, we get the following diagram:
$$
\cxymatrix{
\tilde X\ar[r]^{\tilde\pi}\ar[d]_{\tilde\psi}&X\ar[d]^\psi\\
Y\ar[r]_\pi&Y/\Gamma
}
$$
Since $\tilde\pi$ is finite, it is in fact surjective. Because
$\Gamma$ does not contain reflections, the morphism $\pi$ is \'etale
outside a subset $Y_s\subseteq Y$ of codimension $\ge2$.  Then
$\tilde\pi$ is \'etale outside $X_s:=\tilde\psi^{-1}(Y_s)$. Since
$\psi$ is equidimensional and dominant the same holds for
$\tilde\psi$. This implies that the codimension of $\tilde X_s$ is
$\tilde X$ is at least $2$. Since $X$ is smooth, the Zariski\_Nagata
theorem (stating that the ramification locus is pure of codimension
one, see \cite{SGA2}, Exp.~X, Th.~3.4) implies that $\tilde\pi$ is
(everywhere) \'etale and therefore a covering. The affine space $X$
does not possess non\_trivial coverings. We conclude that
$\tilde\pi$ is an isomorphism. From this, we get a
diagonal morphism $X\pfeil Y$ and therefore
$\CC[Y]^\Gamma\subseteq\CC[Y]\subseteq\CC[X]$. Since $\CC[Y/\Gamma]$
is integrally closed in $\CC[X]$ we conclude $\Gamma=1$.\qed

\Remark: One can replace the Zariski\_Nagata theorem by the following
(well\_known) topological argument: let
$X_s:=\psi^{-1}(Y_s/\Gamma)$. Then $\tilde X_r:=\tilde X\setminus
\tilde X_s \pfeil X_r:=X\setminus X_s$ is a finite \'etale
morphism. The real codimension of $X_s$ in $X$ is at least $4$ and $X$
is $1$\_connected. Therefore, $X_r$ is still $1$\_connected. Because
$\tilde X_r$ is at least connected we conclude that $\tilde X_r\pfeil
X_r$ is an isomorphism. Hence, $\tilde\pi$ is birational and therefore
an isomorphism.

\medskip

Now we obtain that $\psi$ is not only equidimensional but even flat:

\Corollary faithflat. The morphisms $\psi:V\pfeil\fa^*/W_V$ and $\psi\mod
G:V\mod G\pfeil\fa^*/W_V$ are faithfully flat.

\Proof: Since $\psi$ is equidimensional the same holds for $\psi\mod
G$. Since the source of $\psi$ and $\psi\mod G$ is Cohen\_Macaulay
(this is clear for $V$ and follows from the Hochster\_Roberts theorem
for $V\mod G$) and the target is smooth we conclude that both
morphisms are flat (\cite{EGA}~\S15.4.2). This shows in particular
that $\psi(V)$ is open.  On the other hand, this is a homogeneous
subset of $\fa^*/W_V$ which contains $0$. Thus,
$\psi(V)=\fa^*/W_V$.\qed

In algebraic terms, we get:

\Corollary freeness. The rings $\CC[V]$ and $\CC[V]^G$ are free as
$\CC[\fa^*/W_V]$\_modules.

\Proof: For positively graded algebras, freeness and flatness are the same
(\cite{Bour} Ch.~2, \S11, no.~4, Prop.~7).\qed

\beginsection genericstructure. The generic structure

In this section, we refine \cite{generic2} to give a very precise
description of the generic structure of a symplectic representation.
We start with a rather general remark. Recall that an open subset $U$
of an affine $G$\_variety $X$ is called {\it saturated} if it is the
preimage of an open subset of $X\mod G$. This is equivalent to saying
that every closed orbit of $U$ is closed in $X$.

\resetitem

\Proposition selfdual. Let $V$ be a selfdual representation of $G$. Then:
\Item3 For every $G$\_stable divisor $D\subset V$ there is a
non\_zero invariant function $f\in\CC[V]^G$ vanishing on $D$.
\Item4 Every invariant rational function is the quotient of two
regular invariants.
\Item5 Let $F$ be a general fiber of $\pi:V\pfeil V\mod G$. Then $F$
contains an open orbit $F_0$ with $\|codim|_FF\setminus F_0\ge2$.
\Item{33} Let $U\subseteq V$ be a $G$\_stable non\_empty open subset
of $V$. Assume moreover that $U$ is affine. Then $U$ contains a
saturated non\_empty open subset.

\Proof: \cite{I3} The divisor $D$ is the zero set of a function
$f\in\CC[V]$. Moreover, $f$ is unique up to a scalar which implies that
$gf=\chi(g)f$ for all $g\in G$ and a character $\chi$ of $G$. Since
$V$ is selfdual also $\CC[V]_{\le d}=\oplus_{i=0}^dS^iV^*$ is
selfdual. Thus there is a semiinvariant $f^*\in\CC[V]$ with
$gf^*=\chi(g)^{-1}f^*$ for all $g\in G$. Hence, $F=ff^*$ is an
invariant vanishing on $D$.

\cite{I4} It is well known that every invariant rational function $f$ is
the quotient of two regular semiinvariants $f=p/q$ where $p$ and $q$
transform with the same character $\chi$. As above, there is a
semiinvariant $q^*$ with character $\chi^{-1}$. Thus, $f=pq^*/qq^*$ is
a ratio of two invariants.

\cite{I5} Suppose the assertion is false. Then $F$ contains a
$G$\_stable divisor. Using, e.g., the slice theorem this implies that
$V$ contains a $G$\_stable divisor $D$ which meets the general fiber of
$V\pfeil V\mod G$. But then no non\_zero invariant would vanish on
$D$.

\cite{I22} Since $U$ is affine, its complement $D$ in $V$ is a
$G$\_stable divisor. Let $Z:=\overline{\pi(D)}$. Then $Z\ne V\mod G$
by \cite{I3}. Hence $\pi^{-1}(V\mod G - Z)$ has the required
properties.\qed

For symplectic representations we can say much more:

\Theorem generic3. Let $V$ be a symplectic $G$\_representation. Then
there is a commutative diagram
$$29
\cxymatrix{
(G\times^LF)\times U\ar[r]\ar[d]&V_0\ar[d]\inj[r]&V\ar[d]^{\pi_V}\\
U\ar[r]^\alpha\ar[d]&U_0\inj[r]\ar[d]&V\mod G\ar[dl]^{\psi\mod G}\\
\fa^*\ar[r]&\fa^*/W_V
}
$$
where all squares are Cartesian and where $\alpha$ is finite and
\'etale. Here, $L$, $F$, $\fa^*$, and $U$ are as in
\cite{generic2}. Moreover, $U_0$ is an open dense subset of $V\mod G$
and $V_0$ is its preimage in $V$.

\Proof: First, I claim that $\tilde V:=V\times_{\fa^*/W_V}\fa^*$ is
irreducible. In fact, since $\fa^*/W_V$ is an affine
space, $\tilde V$ is a complete intersection. Hence all irreducible
components of $\tilde V$ have the same dimension as $V$ and they all
map dominantly to $V$. It follows that $W_V$ acts transitively on the
set of irreducible components and therefore each component maps
dominantly to $\fa^*$, as well. But the generic fibers of $\tilde
V\pfeil\fa^*$ are the same as the ones of $V\pfeil\fa^*/W_V$ and
therefore irreducible which shows that $\tilde V$ is irreducible.

Let $X:=G\times^LF$. Then \cite{generic2} says that $X\times U$ is an
open subset of $\tilde V$. Clearly, we may shrink $U$ to be
affine. Then the complement $\tilde D$ of $X\times U$ in $\tilde V$ is
a divisor. Since $\tilde V\pfeil V$ is finite, the image $D$ of
$\tilde D$ in $V$ is a divisor as well. From \cite{selfdual}\cite{I33}
we get a non\_empty affine saturated open subset $V_0\subseteq V$ with
$V_0\cap D=\leer$. The preimage of $V_0$ in $\tilde V$ is affine and
is contained in $X\times U$. Hence its complement in $X\times U$ is a
divisor and therefore of the form $X\times U_0$ with $U_0\subseteq U$
open. Replace $U$ by $U_0$. Then the image of $X\times U$ in $V$ is
$V_0$ and the morphism $X\times U\pfeil U$ is finite. Put
$U_0:=V_0\mod G$ and let $\alpha:U\pfeil U_0$ the morphism on
$G$\_quotients. This finishes the construction of diagram \cite{E29}.

The upper right square is Cartesian since $V_0$ is saturated. Let
$\tilde V_0=V_0\times_{\fa^*/W_V}\fa^*$. As an open subset of $\tilde
V$ it is irreducible, as well. Moreover, the morphism $X\times
U\pfeil\tilde V_0$ is both finite  and an open embedding. We conclude
that it is an isomorphism. Thus, the big left square is
Cartesian. Going over to $G$\_invariants we conclude that the lower
left square is Cartesian. This implies, by general nonsense, that the
upper left square is Cartesian.

Finally, by shrinking, if necessary, $U_0$ we may assume that the
image of $U_0$ in $\fa^*/W_V$ lies in the part over which the map
$\fa^*\pfeil\fa^*/W_V$ is unramified. This makes $\alpha$ \'etale.\qed

\Corollary. In the \'etale topology there is a non\_empty open subset
of $V\mod G$ over which the quotient morphism $V\pfeil V\mod G$ is a
trivial fiber bundle with fiber $G\times^LF$.

\Remark: It is follows from Luna's slice theorem that every quotient
map is generically a fiber bundle in the \'etale topology but in
general the trivializing \'etale map cannot be controlled. Therefore,
it might be surprising that for symplectic representations a Galois
covering whose group is a subquotient of the Weyl group suffices.

\Corollary. Let $H$ the principal isotropy group of $V$ (i.e., the
isotropy group of a generic closed orbit). Then there is a Levi
subgroup $L$ with $(L,L)\subseteq H\subseteq L$ and $L/H=A$.

In the next statement, we denote the isotropy group of a non\_zero
vector in $\CC^{2n}$ inside $Sp_{2n}(\CC)$ by $Sp_{2n-1}(\CC)$.

\Corollary. Let $H$ be the generic isotropy group of $V$. Then
$$
H\cong H_0\times Sp_{2m_1-1}(\CC)\times\ldots\times Sp_{2m_s-1}(\CC)
$$
with $H_0$ reductive.

\Remark: A non\_reductive generic isotropy group is a quite
exceptional phenomenon. Clearly $G=Sp_{2m}(\CC)$ acting on
$V=\CC^{2m}$, $m\ge1$ is an example. The same holds more generally for
$G=Sp_{2m}(\CC)\times SO_{2n-1}(\CC)$ acting on
$V=\CC^{2m}\otimes\CC^{2n-1}$ with $2m>2n-1\ge1$. More examples can
be found in \cite{mfclass}.

\beginsection Applications. Symplectic reductions

\Definition: Let $V$ be a symplectic vector space. A fiber of
$\psi\mod G:V\mod G\pfeil\fa^*/W_V$ is called a {\it symplectic
reduction} of $V$.

\Theorem sympred. The symplectic reductions of $V$ form a flat
family parametrized by an affine space. Moreover, for a generic
symplectic reduction there exists a birational Poisson morphism to
$\CC^{2c}$ with its standard symplectic structure.

\Proof: The first part is just a reformulation of parts of
\cite{mainresult} and \cite{faithflat}. The second part follows from
\cite{generic3}.\qed

\Definition: We call $\|rk|_sV:=\|dim|\fa^*$ the {\it (symplectic) rank}
and $c_s(V):=\half(\|dim|V\mod G-\|dim|\fa^*)$ the {\it (symplectic)
complexity} of $V$.

\medskip \noindent Thus the symplectic reductions are of dimension
$2c_s(V)$ and they are parametrized by an affine space of dimension
$\|rk|_sV$.

There already exists a notion of rank and complexity for an
arbitrary $G$\_variety $X$: the complexity $c(X)$ is the transcendence
degree of $k(X)^B$ over $k$ and the rank $\|rk|X$ is the transcendence
degree of $k(X)^U$ minus the complexity of $X$. This previous notion
is related to our present definition:

\Proposition twonotion. For a finite dimensional (non\_symplectic)
representation $X$ of $G$ consider the symplectic representation
$V=X\oplus X^*$. Then $\|rk|_sV=\|rk|X$ and $c_s(V)=c(X)$.

\Proof: Observe that $V=X\oplus X^*$ is just the cotangent bundle
$T^*_X$ of $X$. Then we recognize the assertion as a special case of
\cite{WuM} Thm.~7.1.\qed

\Theorem Poissoncenter. If we regard $\CC[\fa^*]^{W_V}$ as a subalgebra of
$\CC[V]^G$ then it is precisely the Poisson center. In particular,
the Poisson center of $\CC[V]^G$ is a polynomial ring over which it is
a free module.

\Proof: Let $f\in\CC[\fa^*/W_V]$. Then
$f$ is finite over $\CC[\fa^*/W_G]$, i.e., satisfies a monic equation
$p(f)=0$ with coefficients in $\CC[\fa^*/W_G]$. Assume that $p$ is of
minimal degree. Since $\CC[\fa^*/W_G]$ is in the Poisson center of
$\CC[V]^G$ (see section~\cite{momentmap}) we infer for any $g\in\CC[V]^G$:
$$
0=\{p(f),g\}=p'(f)\{f,g\}.
$$
This shows $\{f,g\}=0$ since, by minimality $p'(f)\ne0$. Hence $f$ is
in the Poisson center of $\CC[V]^G$.

Conversely, let $f$ be in the Poisson center of $\CC[V]^G$. Then the
restriction of $f$ to a generic symplectic reduction $R$ is in the
center of $\CC[R]$. According to \cite{sympred}, $R$ has a dense open
subset $R^\circ$ which is a symplectic variety. Since $f$
Poisson\_commutes with all of $\CC[R]$, the same holds for the
localization $\CC(R)$ and therefore $\CC[R^\circ]$. This implies that
$f$ is constant on $R$ and therefore
$f\in\CC[V]^G\cap\CC(\fa^*/W_V)=\CC[\fa^*/W_V]$.\qed

We end this section with a permanence property. It is very useful in
computing complexity and rank of a symplectic representation. See
\cite{mfclass} for examples.

\Theorem. Let $V$ be a symplectic $G$\_representation and let $(S,M)$
be as in \S\cite{local}. Then $\|rk|_s(S,M)=\|rk|_s(V,G)$ and
$c_s(S,M)=c_s(V,G)$. Moreover, there is an $w\in W_G$ with
$\fa^*_V=w\fa^*_S$ and $L_V=wL_Sw^{-1}$.

\Proof: We already noted that $\fa^*$ is, up to conjugation by $W_G$,
determined by the image of $V$ in $\ft^*/W_G$. Thus $\fa^*_V=w\fa^*_S$
follows from \cite{LSThm}. This immediately implies the statement on
ranks. The assertion on $L$ follows from the fact
that $L$ is the point\_wise stabilizer of $\fa^*$. Finally, the
equality of complexities follows, for example from \cite{E33} and the
fact that $\|dim|U=\|rk|_s+2c_s$.\qed

\noindent 

\beginsection multiplicityfree. Multiplicity free symplectic
representations

Of particular interest is the case of complexity zero.

\Definition: A symplectic representation $V$ is called {\it
multiplicity free} if $c_s(V)=0$.

\medskip

\noindent There are many equivalent characterizations of multiplicity
free symplectic representations. Some of them are summarized
below. Observe that \cite{I11} is a precise version of the assertion
``Almost all $G$\_invariants on $V$ are pull-backs of invariants on
$\fg^*$ via the moment map.'' For \cite{I21} recall that a subspace
$U$ of a symplectic vector space is {\it coisotropic} if
$U^\perp\subseteq U$. A smooth subvariety is coisotropic if each of
its tangent spaces is coisotropic.

\resetitem

\Proposition mfPoisson. Let $V$ be a symplectic representation. Then the
following statements are equivalent:
\Item{9} $V$ is multiplicity free.
\Item{23} The morphism $V\mod G\pfeil\fa^*/W_V$ is an isomorphism.
\Item{20} All symplectic reductions are points (in the
scheme sense).
\Item{22} The morphism $V\mod G\pfeil\ft^*/W_G$ is unramified (i.e.,
injective on tangent spaces) on a
dense open subset.
\Item{11} The morphism $V\mod G\pfeil\ft^*/W_G$ is finite.
\Item{12} The Poisson algebra $\CC[V]^G$ is commutative.
\Item{21}The generic $G$\_orbit of $V$ is coisotropic.

\Proof: The morphism $\psi:V\mod G\pfeil\fa^*/W_V$ is faithfully flat of
relative dimension $2c_s(V)$. Thus, \cite{I9} implies \cite{I23} which
implies in turn \cite{I20}, \cite{I22}, and \cite{I11}. Conversely,
each of these conditions imply $c_s(V)=0$ and therefore \cite{I9}.

The equivalence of \cite{I23} and \cite{I12} is a special case of
\cite{Poissoncenter}. Finally, for \cite{I21} consider
$Z:=\overline{m(V)}\subseteq\fg^*$. Then we have a commutative diagram
$$
\cxymatrix{
V\ar[r]^m\ar[d]^{\pi_V}&Z\ \inj[r]\ar[d]^{\pi_Z}&\fg^*\ar[d]^{\pi_{\fg^*}}
\\
V\mod G\ar[r]^{m'}&\fa^*/W_G\inj[r]&\ft^*/W_G
}
$$
where all vertical arrows are categorical quotients.  Choose a generic
point $v\in V$ and denote its image in $Z$, $V\mod G$, and $\fa^*/W_G$
by $z$, $\vq$, and $\zq$, respectively. Let $V_\zq$, $Z_\zq$ etc.
denote the fiber over $\zq$. Then $V_\zq=m^{-1}(Z_\zq)$. The fibers of
$\pi_Z$ contain only finitely many $G$\_orbits (it inherits this
property from $\pi_{\fg^*}$). This implies that $Gz$ is open in
$Z_\zq$. Thus, the preimage of $Gz$ in $V_\zq$ is open as well. This
preimage is isomorphic to the fiber product $G\times^{G_z}V_z$ where
$V_z=m^{-1}(z)$. The tangent space of $V_z$ in $v$ is $(\fg v)^\perp$
(see \cite{E19}).  Thus $T_vV_\zq=\fg v+(\fg v)^\perp$. On the other
hand, $Gv$ is dense in $V_\vq=\pi_V^{-1}(\vq)$
(\cite{selfdual}\cite{I5}). This (and the fact that $V_\zq$ is smooth
in $v$) shows $T_vV_\vq=\fg v$.  Thus, we get
$$
T_\vq(V\mod G)_\zq=(\fg v+(\fg v)^\perp)/\fg v.
$$
This space is zero if and only if $(\fg v)^\perp\subseteq\fg v$,
i.e., $Gv$ is coisotropic. This shows the equivalence of \cite{I22} and
\cite{I21}.\qed

\noindent The notion ``multiplicity free'' is justified by:

\Proposition. For a finite dimensional (non\_symplectic)
representation $X$ of $G$ consider the symplectic representation
$V=X\oplus X^*$. Then $V$ is multiplicity free (as symplectic
representation) if and only if $\CC[X]$ is a multiplicity free
$G$\_module.

\Proof: Using \cite{twonotion}, this is the main statement of
\cite{VK} (see also \cite{Mon} Thm.~3.1).\qed

Associated to a symplectic vector space $V$ is another
object, namely the Weyl algebra $\cW(V)$. By definition, it is an
associative unital algebra which is generated by $V$ with the
relations
$$
uv-vu=\omega(u,v)\quad\hbox{for all $u,v\in V$}.
$$
If $G$ acts on $V$ then it will also act on $\cW(V)$ by way of
automorphisms. On the Lie algebra level, this action is inner: there
is a Lie algebra homomorphism $\rho:\fg\pfeil\cW(V)$ such that
$$
\xi a=[\rho(\xi),a]\quad\hbox{for all $\xi\in\fg$ and $a\in\cW(V)$.}
$$
(It suffices to show this for $\fg=\fs\fp(V)$. In that case see
\cite{Ho} Thm.~5). This means in particular, that the algebra of
invariants $\cW(V)^G$ can be interpreted as centralizer of
$\rho(\fg)$.

The map $\rho$ induces an algebra homomorphism
$\cU(\fg)\pfeil\cW(V)$. With $\cZ(\fg)=\cU(\fg)^G$, the center of
$\cU(\fg)$, we get a homomorphism $\cZ(\fg)\pfeil\cW(V)^G$.

The connection with Poisson algebras is as follows: both algebras
$\cU(\fg)$ and $\cW(V)$ come with natural filtrations by placing the
generators in degree one. For a filtered vector space $F_{\le\bullet}$
let $\|gr|F:=\oplus_{n\in\ZZ}F_{\le n}/F_{<n}$ be the associated
graded space. Then $\|gr|\cU(\fg)=S^*(\fg)=\CC[\fg^*]$ while
$\|gr|\cW(V)=S^*(V)=\CC[V^*]=\CC[V]$. Moreover, the commutator on a
filtered algebra induces a Poisson structure on the associated graded
algebra. In our case, we arrive exactly at the usual Poisson
structures on $\CC[\fg^*]$ and $\CC[V]$. Since $G$ is linearly
reductive we have
$\|gr|\cZ(\fg)=\|gr|\cU(\fg)^G=\CC[\fg^*]^G=\CC[\fg^*\mod
G]$. Similarly, $\|gr|\cW(V)^G=\CC[V]^G$. The homomorphism
$\cZ(\fg)\pfeil\cW(V)^G$ gives on the associated graded level a
homomorphism $\CC[\fg^*\mod G]\pfeil\CC[V]^G$ which corresponds
precisely to the morphism $V\pfeil\fg^*\mod G$ induces by the moment
map.

Now we have further characterizations in terms of the Weyl
algebra. Here statement \cite{I14} is a precise version of the
assertion ``Almost all elements of $\cW(V)$ commuting with $\rho(\fg)$
come from $\cZ(\fg)$''.

\Proposition Weylalgebra.  Let $V$ be a symplectic
representation. Then the following statements are equivalent:
\Item{13} $V$ is multiplicity free.
\Item{14} $\cW(V)^G$ is a finitely generated $\cZ(\fg)$\_module.
\Item{15} The algebra $\cW(V)^G$ is commutative.\Par

\Proof: \cite{I13}$\Rightarrow$\cite{I14}: if $V$ is multiplicity free
then $\CC[V]^G$ is a finitely generated $\CC[\fg^*\mod G]$\_module
(\cite{mfPoisson}~\cite{I11}). This implies the corresponding
statement \cite{I14} on filtered objects.

\cite{I14}$\Rightarrow$\cite{I13}: every $x\in\cW(V)^G$ satisfies a
(monic) equation with coefficients in $\cZ(\fg)$. Looking at highest
degree terms, we see that $\CC[V]^G$ algebraic over $\CC[\fg^*\mod
G]$. Thus, \cite{mfPoisson}~\cite{I22} is satisfied.

\cite{I13}$\Rightarrow$\cite{I15}: we adapt the argument of
\cite{HC}. Suppose $\cW(V)^G$ were not commutative. Let $x\in\cW(V)^G$
be not in the center and consider the derivation $\theta:=\|ad|x$ of
$\cW(V)^G$. Then there is a minimal number $m\in\ZZ$ such that
$\|deg|\theta(y)\le\|deg|y+m$. Thus, $\theta$ induces a non\_zero
derivation $\thetaq$ of $\CC[V]^G$ which is of degree
$m$. Moreover, $\thetaq$ is trivial on $\|gr|\cZ(\fg)=\CC[\fg^*\mod
G]$. By \cite{mfPoisson}~\cite{I11}, that ring is finite in
$\CC[V]^G$. Hence, $\thetaq=0$.

\cite{I15}$\Rightarrow$\cite{I13}: apply
\cite{mfPoisson}~\cite{I12}.\qed

\noindent We finish with a summary of our results specialized to
multiplicity free representations:

\Corollary. Let $V$ be a multiplicity free symplectic representation.
\Item{16} The morphism $V\pfeil\fa^*/W_V$ identifies $\CC[\fa^*]^{W_V}$
with $\CC[V]^G$.
\Item{17} $V$ is a cofree representation, i.e., the invariant ring
$\CC[V]^G$ is a polynomial ring and $\CC[V]$ (or any module of
covariants) is a free $\CC[V]^G$\_module.
\Item{18} The invariant algebra $\cW(V)^G$ is a polynomial ring and
$\cW(V)$ is a free (left or right) $\cW(V)^G$\_module.\Par

\Proof: Part~\cite{I16} is a restatement of
\cite{mfPoisson}\cite{I23}. Part \cite{I17} is just a specialization
of \cite{freeness} to the multiplicity free case. Finally, \cite{I18}
follows from \cite{Weylalgebra}\cite{I15} and its associated graded
version \cite{I17}.\qed

\noindent Because the moment map is quadratic, we get:

\Corollary. Let $V$ be a multiplicity free symplectic
representation. Then the degrees of the generators of $\CC[V]^G$ are
twice the degrees of the generators of $\CC[\fa^*]^{W_V}$. In
particular, they are all even.

\Remark: This last statement is often a convenient way to compute
$W_V$. See \cite{mfclass} where $W_V$ is listed for all cases.

\beginrefs

\B|Abk:Bour|Sig:\\Bo|Au:Bourbaki, N.|Tit:Alg\`ebre, 3. ed%
|Reihe:-|Verlag:Hermann|Ort:Paris|J:1967|xxx:-||

\B|Abk:Bou|Sig:\\Bo|Au:Bourbaki, N.|Tit:Groupes et alg\`ebres de Lie,
Chap. 4, 5 et 6|Reihe:-|Verlag:Hermann|Ort:Paris|J:1968|xxx:-||

\L|Abk:BLV|Sig:BLV|Au:Brion, M.; Luna, D.; Vust, Th.|Tit:Espaces homog\`enes
sph\'eriques|Zs:Invent. Math.|Bd:84|S:617--632|J:1986|xxx:-||

\L|Abk:EGA|Sig:EGA|Au:Dieudonn\'e, J.; Grothendieck, A.|Tit:El\'ements de
g\'eom\'etrie alg\'ebrique IV|Zs:Publ. Math. IHES|Bd:28|S:-|J:1966|xxx:-||

\B|Abk:GS|Sig:GS|Au:Guillemin, V.; Sternberg, S.|Tit:Symplectic
techniques in physics|Reihe:-|Verlag:Cambridge University
Press|Ort:Cambridge|J:1984|xxx:-||

\L|Abk:Ho|Sig:Ho|Au:Howe, R.|Tit:Remarks on classical invariant
theory|Zs:Trans. Amer. Math. Soc.|Bd:313|S:539--570|J:1989|xxx:-||

\L|Abk:WuM|Sig:\\Kn|Au:Knop, F.|Tit:Weylgruppe und Momentabbildung|Zs:Invent.
Math.|Bd:99|S:1--23|J:1990|xxx:-||

\L|Abk:HC|Sig:\\Kn|Au:Knop, F.|Tit:A Harish\_Chandra homomorphism for
reductive group actions|Zs:Annals Math. (2)|Bd:140|S:253--288|J:1994|xxx:-||

\Pr|Abk:Mon|Sig:\\Kn|Au:Knop, F.|Artikel:Some remarks on multiplicity
free spaces|Titel:Proc. NATO Adv. Study Inst. on Representation Theory and
Algebraic Geometry|Hgr:A.~Broer, G.~Sabidussi, eds.|Reihe:Nato ASI
Series C|Bd:514|Verlag:Kluwer|Ort:Dortrecht|S:301--317|J:1998|xxx:-||

\L|Abk:mfclass|Sig:\\Kn|Au:Knop, F.|Tit:Classification of multiplicity
free symplectic representations|Zs:Preprint|Bd:-|S:25
pages|J:2005|xxx:math.SG/0505268||

\L|Abk:LuAO|Sig:Lu|Au:Luna, D.|Tit:Adh\'erences d'orbite et
invariants|Zs:Invent. Math.|Bd:29|S:231--238|J:1975|xxx:-||

\L|Abk:Pan|Sig:Pa|Au:Panyushev, D.|Tit:On orbit spaces of finite and
connected linear groups|Zs:Math. USSR Izv.|Bd:20|S:97-101|J:1983|xxx:-||

\B|Abk:SGA2|Sig:SGA2|Au:Grothendieck, A.|Tit:Cohomologie locale des
faisceaux coh\'erent et th\'eor\`emes de Lefschetz locaux et
globaux|Reihe:-|Verlag:Masson et Cie, North
Holland|Ort:Amsterdam|J:1968|xxx:-||

\L|Abk:VK|Sig:VK|Au:Vinberg, E., Kimelfeld,
B.|Tit:Homogeneous domains on flag manifolds and spherical subsets of
semisimple Lie groups|Zs:Funktsional. Anal. i
Prilozhen.|Bd:12|S:12--19|J:1978|xxx:-||

\endrefs

\bye